\documentclass[]{interact}

\usepackage{amsthm}
\usepackage{amssymb} %解决 /mathbb 的问题
\usepackage{amsmath}
\usepackage{fancyhdr}
\usepackage{hyperref}
\usepackage{mathtools}
\usepackage{physics}
\usepackage{color}

\numberwithin{equation}{section}%按照章节编号

\newtheorem{definition}{Definition}[section]
\newtheorem{theorem}[definition]{Theorem}
\newtheorem{lemma}[definition]{Lemma}
\newtheorem{proposition}[definition]{Proposition}

\newtheorem{remark}{Remark}[section]

\title{Upper capacity entropy and packing entropy of saturated sets for amenable group actions}

\begin{document}

\bibliographystyle{plain}

\author{
	\name{Xiankun Ren\textsuperscript{a,$^\ast$}, Wenda Zhang\textsuperscript{b,$^{\S}$}, and Yiwei Zhang\textsuperscript{c,$^{\ddag}$}
		\thanks{\textsuperscript{$^\ast$}Email:xkren@cqu.edu.cn; \textsuperscript{$^{\S}$} wendazhang951@aliyun.com; \textsuperscript{$^{\ddag}$}yiweizhang@hust.edu.cn 
		\newline 2010 Mathematics Subject Classification.  37B05,37B40,54H20
	\newline \textsuperscript{$^{\S}$} The second author is the corresponding author}
	}	
	\affil{\textsuperscript{a}College of Mathematics and Statistics\\Chongqing University
		\\ Chongqing 401331, China}	
	\affil{\textsuperscript{b}College of Mathematics and Statistics, Chongqing Jiaotong University
		\\ Chongqing 400074, China}
	\affil{\textsuperscript{c}School of Mathematics and Statistics, Center for
		Mathematical Sciences, Hubei Key Laboratory of Engineering Modeling
		and Scientific Computing, Huazhong University of Sciences and
		Technology, Wuhan 430074, China}
}

\maketitle

\begin{abstract}
Let $(X,G)$ be a $G$-action topological system, where $G$ is a countable infinite discrete amenable group and $X$ a compact metric space. In this paper we study the upper capacity entropy and packing entropy for systems with weaker version of specification. We prove that the upper capacity always carries full entropy while there is a variational principle for packing entropy of saturated sets.
\end{abstract}


\section{Introduction}

The study of the size of saturated sets for dynamical systems was initialized by Bowen  and later was developed by Sigmund, Pfister and Sullivan  et. al. This paper focuses on establishing upper capacity and packing entropy formulae of saturated sets for countable amenable group actions with specification or weaker version.

Throughout this paper, $X$ will be a compact metric space with meric $\rho$ and $G$ a discrete countable amenable group. Let $M(X), M(X,G)$ and $E(X,G)$ denote  the sets of all Borel probability measures, $G$-invariant Borel probability measures and ergodic $G$-invariant Borel probability measures induced with the weak$^*$ topology, respectively . When $G=\mathbb{Z}$, we write the system $(X,G)$ as $(X,T)$, when $T:X \rightarrow X$ is a homeomorphism. For finite subsets $\Gamma\subset X$, we denote $\#\Gamma$ the cardinality of $\Gamma$.

Let $\{F_n\}$ be a F{\o}lner sequence. The {\it empirical measure} of $x\in X$ with respect to one finite subset $F\subset G$ is the Borel probability measure
\[
\mathcal{E}_{F}(x) := \frac{1}{|F|}\sum_{s\in F} \delta_{sx},
\]  
where $\delta_y$ is the Dirac measure at point $y$. Let $V(x,\{F_n\})$ be the limit-point set of  $\{\mathcal{E}_{F_n}(x)\}$, which is always a compact connected subset of $M(X,G)$. For each compact connected non-empty subset $K\subset M(X,G)$, denote by $Y_K(\{F_n\}) = \{x\in X \mid V(x,\{F_n\}) = K\}$ which is called {\itshape the saturated set of $K$.} We write $Y_K(\{F_n\})$ as $Y_K$ for simplicity. In particular, when specializing $K$ to be a singleton measure $\mu\in M(X,G),$ the set $Y_{K}$ coincides with the set of generic points of $\mu$ with respect to $\{F_n\}$ and we write $Y_{\{\mu\}}$ as $Y_\mu$.

Bowen introduced a definition of topological entropy of subsets in \cite{bowen1973topological} which is known as Bowen topological entropy and obtained a remarkable {\itshape single-saturated property} 
\[
h_{top}^{B}(G_{\mu}, T) = h_{\mu}(X, T) \text{ when }\mu \text{ is ergodic},
\]
where $h_{\mu}(T,X)$ is the metric entropy of $\mu$. The interesting part of Bowen's formula is that the metric entropy is affine while Bowen topological entropy can be viewed as dynamically analogous to Hausdorff dimension.

After Bowen's work, a series of papers pursued Bowen's formula further into more general saturated set $Y_K$, with $K$ to be a compact connected subset in $M(X,T)$.  Pister and Sullivan \cite{pfister2007topological} pointed out that it is not difficult to provide examples with $h_{\mu}(X,T)>0$ but the generic set of $\mu$ is empty. However, Denker, Grillenberg and Sigmund \cite{denker2006ergodic} proved that $Y_{K}$ is not empty foy systems with hyperbolicity or specification and this result is generalized to non-uniformly hyperbolic systems in \cite{liang2013ergodic}. Pfiter and Sullivan proved a {\itshape saturaed property} for $T$ satisfying {\itshape the $g$-almost product property} and {\itshape uniform separation property}(see the definitions in Section 2) the following holds:
\begin{align}
h_{top}^{B}(Y_K, T) = \inf\{h_{\mu}(X,T)\mid \mu\in K\}. \label{intr1}
\end{align}
The $g$-almost product property is weaker than specification and  uniform separation is weaker than expansiveness. Also the equation \ref{intr1} is generalized into non-uniformly hyperbolic systems in \cite{liang2017variational} and non-uniformly expanding maps in \cite{tian2017topological}.

Besides Bowen topological entropy, {\itshape packing topological entropy} and {\itshape upper capacity topological entropy} are another two important concepts for characterizing the size of non-compact subset $Z\subset X$. Actually, the upper capacity topological entropy is a direct generalization of the classical topological entropy \cite{adler1965topological}. The packing topological entropy was introduced by Feng and Huang \cite{feng2012variational} in a way by resembling packing dimension in dimension theory. 

It is natural to ask if  there are different saturated properties for Bowen topological entropy, packing topological entropy and upper capacity topological entropy. The authors in \cite{hou2020saturated} gave answers to this question for $\mathbb{Z}$-actions. In this paper we are dealing with this question for more general group actions.

In contrast to amenable group $\mathbb{Z}$, a more general countable amenable may have more complicated structure and new phenomena and difficulties may arise.  In 1987, Ornstein and Weiss \cite{ornstein1987entropy} developed the so-called {\itshape quasi-tiling} method, which has been a basic tool in the study of amenable group actions. Many researchers  have made lots of progress in many directions of ergodic theory for amenable group actions. For example Kieffer \cite{kieffer1975generalized} proved the Shannon-McMillan theorem for amenable group actions; Lindenstrauss \cite{lindenstrauss1999pointwise} later established this pointwise result for tempered F{\o}lner sequence with superlogarithmic growth, which can be found in any amenable group. Downarowicz, Huczek and Zhang \cite{downarowicz2019tilings} recently have a new result on the quasi-tiling. Also Zhang \cite{zhang2018topological} used the new quasi-tiling techniques to calculate  topological pressure of generic points.
For more information of amenable group actions, readers may refer \cite{huang2011local,ollagnier2007ergodic,kerr2016ergodic}.

Our main results are the following theorems.
\begin{theorem}\label{mainthm1}
	Suppose the topological dynamical system $(X,G)$ satisfies the $g$-almost product property. Let $K\subset M(X,G)$ be a compact connected nonempty set and $\{F_n\}$ be a F{\o}lner sequence. %with $\frac{|F_n|}{\log n} \rightarrow \infty.$
	Then we have
	\[
	h^{UC}_{top}(Y_{K},\{F_n\}) = h_{top}(X,\{F_n\}) = h_{top}(X,G).
	\]
\end{theorem}

\begin{theorem}\label{mainthm2}
	Suppose the topological dynamical system $(X,G)$ satisfies the $g$-almost product property and has uniform saparation property. Let $K\subset M(X,G)$ be a compact convex nonempty set and $\{F_n\}$ be a F{\o}lner sequence with $\frac{|F_n|}{\log n} \rightarrow \infty.$  Then we have
	\[
	h^{P}_{top}(Y_K,\{F_n\}) = \sup\{h_{\mu}(X,G)\mid \mu\in K\}.
	\]
\end{theorem}

From Theorem \ref{mainthm1} and Theorem \ref{mainthm2}, we notice that the upper capacity topological entropy and the packing topological entropy show quite different behavior from saturated property viewpoint.

The paper is organized as follows. In Section 2, we introduce some notions and lemmas. In Section 3, we prove the Theorem \ref{mainthm1} and Theorem \ref{mainthm2}. In section 4, we give an application of the main theorems. In Appendix, we give the proof of Theorem \ref{approthm}.
\section{Preliminary}

In this section, we will introduce some notions and recall some lemmas.

Let  $F(G)$ be the collection of all finite subsets of $G$. Let $K\in F(G)$. The $K$ boundary of $F\in G$ is defined by 
\[
\partial_{K}(F) := \{c\in G \mid Kc\cap F \neq\emptyset, Kc\cap(G\setminus F) \neq \emptyset\}.
\]
A set $F\in F(G)$ is called $(K,\delta)-invariant$ if $\frac{|\partial_{K}(F)|}{|F|} < \delta.$ A sequence $\{F_n\}\subset F(G)$ is called a {\it F{\o}lner sequence} if for any $s\in G$,
\[
\lim\limits_{n\rightarrow\infty}\frac{|gF_n\triangle F_n|}{|F_n|} = 0.
\]
And $\{F_n\}$ is a F{\o}lner sequence if and only if for any $K\in F(G)$ and $\delta>0$, there exists $N>0$ such that when $n>N$, $F_n$ is $(K,\delta)$-invariant.

A F{\o}lner sequence $\{F_n\}_{n=1}^{\infty}$ is said to be {\itshape tempered} if there exists $C>0$ which is independent of $n$ such that
\[|\bigcup_{k<n}F_{k}^{-1}F_n| \leq C|F_n|.\]
A F{\o}lner sequence is said to have the {\it superlogarithmic growth}  if
\begin{align}
\lim\limits_{n\rightarrow\infty}\frac{|F_n|}{\log n} = \infty. \label{increaing}
\end{align}
For $F\in F(G)$, let $\rho_F$ be the metric defined by 
\[
\rho_{F}(x,y) := \max_{s\in F} \rho(sx,sy).
\]
For $\varepsilon>0, x\in X$ and $F\in F(G)$, we denote by 
\[
B_{F}(x,\varepsilon) := \{y\in X\mid \rho_{F}(x,y)< \varepsilon\}
\]
and
\[
\overline{B}_{F}(x,\varepsilon) := \{y\in X\mid \rho_{F}(x,y)\leq \varepsilon\}
\]
which are the open and closed $F$- {\itshape Bowen balls} with center $x$ and radius $\varepsilon$, respectively.

For $\varphi \in C(X,\mathbb{R})$ and $\mu\in M(X)$, we denote $\langle\varphi,\mu\rangle = \int_{X}\varphi \dd\mu$. There is a countable and separating set of continuous functions $\{\varphi_1,\varphi_2,\cdots\}$ with $0\leq \varphi_i\leq 1, i=1,2,\cdots$ such that 
\begin{align}
D(\mu,\nu) := \sum_{i=1}^{\infty}\frac{|\langle\varphi_i,\mu -\nu\rangle|}{|2^i|}	\label{distance1}
\end{align}
which defines a compatible metric for the weak$^*$-topology on $M(X)$. For $\xi>0$, we denote a ball in $M(X)$ by
\[
B(\mu,\xi) := \{\nu\in M(X)\mid D(\mu,\nu) < \xi\}.
\]

We mention that in the rest of this paper, we will always use the equivalent metric 
\begin{align}
\rho(x,y) := D(\delta_x,\delta_y) \label{metric}
\end{align}
as the metric on $X$, where $\delta_x$ is the Dirac mass at $x\in X.$

\subsection{Entropies for subsets} 

In this part, we introduce three important types of topological entropies for subsets. They are the Bowen topological entropy, packing topological entropy and the upper capacity entropy. The former two are dimension-like concepts and the third one is a topological concept.

\begin{definition}[Bowen topological entropy]

Let $Z\subset X$ and $\{F_n\}$ be a F{\o}lner sequence. For $\varepsilon>0$ and $N\in \mathbb{N},$ denote $\mathcal{C}_{N}(Z,\varepsilon,\{F_n\})$  the collection of all finite or countable cover $\mathcal{C}=\{B_{F_{n_i}}(x,\varepsilon)\}$ of $Z$ with $n_i\geq N.$ For $s>0$, denote
\[
M(Z,\varepsilon,N,s,\{F_n\}) := \inf\limits_{\mathcal{C}\in \mathcal{C}_{N}(Z,\varepsilon,\{F_n\})} \sum_{B_{F_m}(x,\varepsilon)\in \mathcal{C}}e^{-s|F_m|}.
\]
The value $M(Z,\varepsilon,N,s,\{F_n\})$ does not decrease as $N$ increases, hence the following limit exists 
\[
M(Z,\varepsilon,s,\{F_n\}) = \lim\limits_{N\rightarrow \infty}M(Z,\varepsilon,N,s,\{F_n\}).
\]
There exists a critical value of $s$ such that $M(Z,\varepsilon,s,\{F_n\})$ jumps from $+\infty$ to $0.$
Let 

\begin{align*}
h_{top}^{B}(Z,\varepsilon,\{F_n\}) &= \inf \{s\mid M(Z,\varepsilon,s,\{F_n\}) = 0 \}\\
& = \sup \{s\mid M(Z,\varepsilon,s,\{F_n\}) = \infty \}.
\end{align*}
Clearly, $h_{top}^{B}(Z,\varepsilon,\{F_n\})$ dose not decrease as $\varepsilon$ decreases, hence the following limit exists

\[h_{top}^{B}(Z,\{F_n\}) = \lim\limits_{\varepsilon \rightarrow 0}h_{top}^{B}(Z,\varepsilon,\{F_n\}),
\]
and we call it Bowen topological entropy of $Z$ with respect to the sequence $\{F_n\}.$
\end{definition}

The packing topological entropy for amenable group actions was introduced by Dou, Zheng and Zhou \cite{dou2020packing}.

\begin{definition}[Packing topological entropy]
	For $Z\subset X, s\geq 0, N\in \mathbb{N}$ and $\varepsilon>0$, define
	\[
	P(Z,\varepsilon,N,s,\{F_n\}) = sup \sum_{i}\exp(-s|F_{n_i}|),
	\]
	where the supremum is taken over all finite or countable pairwise disjoint families $\{\overline{B}_{F_{n_i}}(x_i,\varepsilon)\}$ with $x_i \in Z, n_i\geq N$ for all $i$. The quantity $P(Z,\varepsilon,N,s,\{F_n\})$ dose not increase as $N$ increases, hence the following limit exists:
	\[
	P(Z,\varepsilon,s,\{F_n\}) = \lim\limits_{N\rightarrow \infty}P(Z,\varepsilon,N,s,\{F_n\}).
	\]
Define 
	\[
	\mathcal{P}(Z,\varepsilon,s,\{F_n\}) = \inf\{\sum_{i=1}^{\infty}P(Z_i,\varepsilon,s,\{F_n\}):\cup_{i=1}^{\infty}Z_i \supset Z\}.
	\]
There exists a critical value of the papameter $s$, which we will denote by $h_{top}^{P}(Z,\varepsilon,\{F_n\})$, where $\mathcal{P}(Z,\varepsilon,s,\{F_n\})$ jumps from $+\infty$ to $0$, i.e. 
	
	\[
	\mathcal{P}(Z,\varepsilon,s,\{F_n\}) = 
	\begin{cases}
	0& s>h_{top}^{P}(Z,\varepsilon,\{F_n\})\\
	+\infty& s<h_{top}^{P}(Z,\varepsilon,\{F_n\}).
	\end{cases}
	\]
It is not hard to see that $h_{top}^{P}(Z,\varepsilon,\{F_n\})$ increases when $\varepsilon$  decreases. We call 
	\[
	h_{top}^{P}(Z,\{F_n\}) := \lim\limits_{\varepsilon \rightarrow 0}h_{top}^{P}(Z,\varepsilon,\{F_n\})
	\]
	the {\itshape amenable packing topological entropy}(packing entropy for short) of $Z$ (with respect to the F{\o}lner sequence $\{F_n\}$).
\end{definition}

Packing entropy has the following properties.
\begin{proposition}\cite[Proposition 2.1]{dou2020packing} \label{packingproperty}
\begin{enumerate}
	\item If $Z^{\prime}\subset Z \subset X$, then
	\[
	h_{top}^{P}(Z^{\prime},\{F_n\}) \leq h_{top}^{P}(Z,\{F_n\});
	\] 
	\item If $Z\subset \bigcup_{i=1}^{\infty}Z_i$, then 
	\[
	h_{top}^{P}(Z,\{F_n\}) \leq \sup_{i\geq 1} h_{top}^{P}(Z_i,\{F_n\});
	\]
	
	\item If $\{F_{n_k}\}$ is a subsequence of $\{F_n\}$, then
	\[
	h_{top}^{P}(Z,\{F_{n_k}\}) \leq h_{top}^{P}(Z,\{F_n\}).
	\]
\end{enumerate}
\end{proposition}

Let $Z\subset X$ be a non-empty set and $F\subset G$ a non-empty finite set. For $\varepsilon>0$, a set $E\subset Z$ is called  $(F,\varepsilon)$-separated if $x\neq y \in E$ implies $\rho_{F}(x,y)>\varepsilon;$ a set $E$  is said to $(F,\varepsilon)$ span $Z$ if for any $y\in Z$, there exists $x\in E$ such that $\rho_{F}(x,y) \leq \varepsilon$. Let $s_{F}(Z,\varepsilon)$ denote the largest cardinality of $(F,\varepsilon)$-separated sets for $Z$ and $r_{F}(Z,\varepsilon)$ the smallest cardinality of any $(F,\varepsilon)$-spanning set for $Z$. 
\begin{definition}[Upper capacity topological entropy]
	Let $\{F_n\}$ be a F{\o}lner sequence. Then  upper capacity topological entropy of $Z$ is defined by
	\[
	h_{top}^{UC}(Z,\{F_n\}) = \lim\limits_{\varepsilon \rightarrow 0}\limsup_{n\rightarrow \infty}\frac{1}{|F_n|}\log s_{F_n}(Z,\varepsilon) = \lim\limits_{\varepsilon \rightarrow 0}\limsup_{n\rightarrow \infty}\frac{1}{|F_n|}\log r_{F_n}(Z,\varepsilon).
	\]
\end{definition}

The following propositions show order relations of $h_{top}^{B}, h_{top}^{P}$ and $h^{UC}_{top}$.

\begin{proposition}\cite[Proposition 2.2]{dou2020packing}
	For any $Z\subset X$ and any F{\o}lner sequence $\{F_n\}$, 
	\[
	h_{top}^{B}(Z,\{F_n\}) \leq h_{top}^{P}(Z,\{F_n\}).
	\]
\end{proposition}

\begin{proposition}\label{packinguc}\cite[Propersition 2.4]{dou2020packing}
	If the F{\o}lner sequence $\{F_n\}$ satisfies $\frac{|F_n|}{\log n}\rightarrow \infty$, then for any subset $Z\subset X$,
	\[
	h_{top}^{P}(Z,\{F_n\}) \leq h_{top}^{UC}(Z,\{F_n\}).
	\]
\end{proposition}

%\begin{proposition}
%	Let $\{F_n\}$ be a F{\o}lner sequence satisfying the increasing condition and $Z\subset X$ a subset of $X$. Then
%	\[
%	h_{top}^{P}(Z,\{F_n\}) = \inf\{\sup_{i\geq 1}h_{top}^{UC}(Z_i,\{F_n\}) : Z=\cup_{i=1}^{\infty}Z_i\}
%	\]
%\end{proposition}

From measure theoretical viewpoint, Brin and Katok \cite{brin1983onlocal} introduced the {\itshape local measure-theoretical lower and upper entropies} of $\mu\in M(X)$. We give the definitions with related to one F{\o}lner sequence $\{F_n\}$.

\begin{definition}\label{localentropy}
	For $\mu \in M(X), x\in X, \varepsilon>0$ and a F{\o}lner sequence $\{F_n\}$, let
	\[
	\overline{h}^{loc}_{\mu}(x,\varepsilon, \{F_n\}) = \limsup_{n\rightarrow \infty} -\frac{1}{|F_n|} \log \mu(B_{F_n}(x,\varepsilon)),
	\]
	
	and
	\[
	\underline{h}^{loc}_{\mu}(x,\varepsilon, \{F_n\}) = \liminf_{n\rightarrow \infty} -\frac{1}{|F_n|} \log \mu(B_{F_n}(x,\varepsilon)).
	\]

For  $Z\in \mathcal{B}(X)$ (the Borel $\sigma$-algebra of $X$), denote by 
	\begin{align*}
	\overline{h}^{loc}_{\mu}(Z,\{F_n\}) = \int_{Z}\lim\limits_{\varepsilon \rightarrow 0} \overline{h}^{loc}_{\mu}(x,\varepsilon, \{F_n\}) \dd \mu
	\end{align*}
	and 
	\begin{align*}
	\underline{h}^{loc}_{\mu}(Z,\{F_n\}) = \int_{Z}\lim\limits_{\varepsilon \rightarrow 0}\underline{h}^{loc}_{\mu}(x,\varepsilon, \{F_n\}) \dd \mu,
	\end{align*}
which are called the {\it upper local entropy} and the {\it lower local entropy} of $\mu$ over $Z$, respectively. 
\end{definition}

The following theorem is an amenable version of \cite[Theorem 1.3]{feng2012variational}. It
 shows a variational principle between packing entropy and local measure-theoretical upper entropy.
\begin{theorem}\cite[Theorem 1.3]{dou2020packing}\label{packingvp}
	Let $(X,G)$ be a $G$-action topological dynamical system and $G$ a discrete countable amenable group. Let $\{F_n\}$ be a F{\o}lner sequence with $\frac{|F_n|}{\log n} \rightarrow \infty$. Then for any non-empty Borel subset $Z\subset X$,
	\[
	h_{top}^{P}(Z,\{F_n\}) = \sup\{\overline{h}_{\mu}^{loc}(Z,\{F_n\})\mid \mu\in M(X),\mu(Z)=1\}.
	\] 
\end{theorem}

%\subsection{Quasi-tilling for amenable group}
%In this part, we will recall the quasi-tilling method and a recent result about tillings for amenable groups from \cite{downarowicz2019tilings}.

\subsection{The $g$-almost product property}
The specification property was first introduced by Bowen. The $g$-almost product property was introduced by Pfister and Sullivan in \cite{pfister2007topological} and it turns out to be a weaker concept than specification property. For instance, it is well known that the $g$-almost product property holds for all $\beta$-shifts. However, the specification property is not satisfied for every parameter of $\beta$-shifts. The $g$-almost product property for amenable group actions was introduced by Zhang in \cite{zhang2018topological}.

\begin{definition}\label{almostspecification}
	A map $g: (0,1) \rightarrow (0,1)$ is called a {\itshape mistake-density function} if 
	\[
	\lim\limits_{r\rightarrow 0}g(r) = 0 \text{ and } g(r)\leq g(s),\quad \forall\  0<r<s<1.
	\]
	For $F\in F(G), \varepsilon>0$ and $x\in X$, define the dynamical ball with respect to the mistake-density function $g$ by
	\begin{align}
			B(g; F, x, \varepsilon) := \{y\in X : |\{s\in F\mid \rho(sx,sy) > \varepsilon\}| \leq g(\varepsilon)|F|\}. \label{dynamicalball}
	\end{align}

\end{definition}

Next we will recall the $g-almost$ product property for group actions.

\begin{definition}
	Let $g$ be a mistake-density function. The system $(X,G)$ satisfies the $g$-almost product property if there exists a map $m:(0,1) \rightarrow F(G) \times (0,1)$ such that for any $k\in \mathbb{N}$, any $\varepsilon_1, \varepsilon_2, \cdots, \varepsilon_k$ and any $x_1, x_2,\cdots, x_k \in X$, if $\{F_i\}_{i=1}^{k}$ are pairwise disjoint and $F_i$ is $m(\varepsilon_i)$-invariant, $i=1,2,\cdots,k$, then
	\[
	\bigcap_{i=1}^{k}B(g; F_i, x_i, \varepsilon_i) \neq \emptyset. 
	\]

\end{definition}

%The following lemma can be easily checked from the  distance \ref{metric} on $X$.
%\begin{lemma}\label{eslemma0}
%	Let $\delta>0, F\in F(G)$ and $x\in X.$ Suppose $F^{\prime} \subset F$ with $\frac{|F^{\prime}|}{|F|} > 1-\delta$. Then 
%	\[
%	D(\mathcal{E}_{F^{\prime}}(x), \mathcal{E}_{F}(x)) < 2\delta.
%	\]
%\end{lemma}

\begin{lemma}\label{eslma1}
	Assume the system $(X, G)$ has the  $g$-almost product property. Let $x_1,x_2,\cdots,x_k \in X$ and $\varepsilon_1>0,\varepsilon_2>0,\cdots,\varepsilon_k>0$ and $F_1,F_2,\cdots,F_k\in F(G)$ be given satisfying each $F_j$ is $m(\varepsilon_i)$-variant for $j=1,2,\cdots,k$ and $F_{i}\cap F_j = \emptyset$ for $1\leq i\neq j \leq k$. Let $F=\bigcup_{i=1}^{k}F_i$. Assume that $\mathcal{E}_{F_j}(x_j) \in B(\mu_j, \xi_j), \ j=1,2,\cdots,k$. Then for any
	$y \in \bigcap_{j=1}^{k}B(g;F_j,x_j,\varepsilon_j)$ and any probability measure $\alpha,$ 
	\[
	D(\mathcal{E}_{F}(y),\alpha) \leq \sum_{i=1}^{k}\frac{|F_i|}{|F|}\big(D(\mu_j,\alpha) + \xi_j +\varepsilon_j+g(\varepsilon_j)\big).
	\]
\end{lemma}

\begin{proof}
	We have
	\[
	\mathcal{E}_{F}(y) = \sum_{j=1}^{k}\frac{|F_i|}{|F|}\mathcal{E}_{F_i}(y).
	\]
	Because of the distance (\ref{metric}) on $X$, 
	\[
	D(\mathcal{E}_{F_j}(x_j),\mathcal{E}_{F_j}(y)) < g(\varepsilon_j) + \varepsilon_j.
	\]
	
	The result follows from the triangle inequality and the definition of the distance (\ref{metric}):
	\[
	D(\mathcal{E}_{F}(y),\alpha) \leq \sum_{j=1}^{k}\frac{|F_j|}{|F|}D(\mathcal{E}_{F_j}(y),\alpha) \leq \sum_{j=1}^{k}\frac{|F_j|}{|F|}\big(D(\mu_j,\alpha) + g(\varepsilon_j) + \varepsilon_j +\xi_j\big).
	\]
\end{proof}
\subsection{Uniform separation property and entropy dense}

Let $F\in F(G)$ and $\delta>0, \varepsilon>0$. A subset $\Gamma \subset X$ is $(\delta,F,\varepsilon)-separated$ if for $x\neq y \in \Gamma,$
\[
\frac{|\{s\in F\mid \rho(sx,sy)>\varepsilon\}|}{|F|} \geq \delta.
\]
Let $C\subset M(X)$ be a neighborhood of $\mu\in M(X,G).$ Define
\[
X_{F,C} := \{x\in X\mid \mathcal{E}_{F}(x)\in C\},
\]
and 
\[
N(C;F,\varepsilon) := \text{largest cardinality of an } (F,\varepsilon)-separated \text{ subset of } X_{F,C} 
\]
and
\[
N(C;\delta, F,\varepsilon) := \text{largest cardinality of a } (\delta, F,\varepsilon)-separated \text{ subset of } X_{F,C}.
\]

With this convention, we have the following proposition.
\begin{proposition}\cite[Propersition 3.7]{ren2020topological}\label{sepes1}
	
	Let $\{F_n\}$ be a  F{\o}lner sequence and $\mu$ be an ergodic measure. Then for $h^*<h_{\mu}(X,G)$, there exist $\delta^*>0\text{ and }\varepsilon^*>0$ such that for any neighborhood $C$ of $\mu$, there exists $n_{C}^*$ , such that for any $n\geq n_{C}^*$ there exists a $(\delta^*,F_n,\varepsilon^*)$-separated set $\Gamma_n$ of $X_{F_n,C}$ satisfying
	\[
	\# \Gamma_n\geq e^{h^*|F_n|}.
	\]
\end{proposition}

The uniform separation property for $\mathbb{Z}$-actions was introduced by Pfister and Sullivan in \cite{pfister2007topological}.
The concept for amenable group actions was defined in \cite[Section 3]{ren2020topological}.
\begin{definition}
	The  system $(X,G)$ has {\itshape uniform separation property} if the following holds. Let $\{K_n\}$ be a tempered F{\o}lner sequence. For any $\eta >0,$ there exists $\varepsilon^*>0$ and $\delta^*>0$ such that for $\mu\in E(X,G)$ and any neighborhood $C\subset M(X)$ of $\mu$, there exists $n_{C;\mu,\eta}^{*}\in\mathbb{N},$ such that for $n\geq n_{C;\mu,\eta}^{*}$,
	\[
	N(C;\delta^*,K_n,\varepsilon^*) \geq e^{|K_n|(h_{\mu}(X,G) - \eta)}.
	\]
\end{definition}

\noindent{\itshape Remark:} The uniform separation property implies $h_{top}(X,G) < \infty.$ 

At last, we will recall the concept of	entropy dense.

\begin{definition}
	The measure $\nu\in M(X,G)$ is entropy-approachable by ergodic measures if for any neighborhood $C$ of $\nu$ and each $h^*< h_{\nu}(X,G),$ there exists a measure
	$\mu\in E(X,G)\cap C$ such that $h_{\mu}(X,G) > h^*.$ The ergodic measures are  entropy-dense if each $\nu\in M(X,G)$ is entropy-approachable by ergodic measures.
\end{definition}

For amenable group actions, we prove the following result.

\begin{theorem}\label{approthm} Suppose the system $(X,G)$ has  the $g$-almost product property. Then the ergodic measures are entropy dense.
\end{theorem}
\begin{proof}
	We will prove this theorem in the Appendix part.	
\end{proof}

The following proposition is \cite[Proposition 3.4]{ren2020topological}
\begin{proposition}\label{separatinges1}
	Let $(X,G)$ be a dynamical system. Assume the system has the uniform separation property and that the ergodic measures are entropy dense. Let $\{K_n\}$ be a tempered F{\o}lner sequence. For any $\eta >0$, there exist $\delta^*>0$ and $\varepsilon^*>0$ so that for $\mu\in M(X,G)$ and any neighborhood $C\subset M(X)$ of $\mu$, there exists $n^{*}_{C;\mu,\eta}$ such that for $n\geq n^{*}_{C;\mu,\eta}$
	\[
	N(C;\delta^*,K_n,\varepsilon^*) \geq e^{|K_n|(h_{\mu}(X,G) - \eta)}.
	\]

\end{proposition}

The following lemmas will be needed.

\begin{lemma}\cite[Lemma 3.6]{ren2020topological}\label{upperlemma}
	Let  $\mu\in M(X,G)$ and $\{F_n\}$ be a F{\o}lner sequence. For $\varepsilon>0$, let $E_n$ be a sequence of $(F_n, \varepsilon)$-separated subsets. Define
	$$\nu_{n} := \frac{1}{|E_n|}\sum_{x\in E_n}\mathcal{E}_{F_n}(x).$$
	
	Assume $\nu_n \rightarrow \mu.$ Then
	$$\limsup_{n\rightarrow \infty.}\frac{1}{|F_n|}\log|E_n| \leq h_{\mu}(X,G).$$
\end{lemma}

\begin{lemma}\cite[Lemma 2.6]{zhang2018topological}\label{eslemma1}
	
	Let $(X,G)$ be a dynamical system. Let $\mu\in M(X,G), \delta^{*}>0, \varepsilon^{*} >0,\xi>0.$ Let $0< \delta < \min\{\frac{1}{2},\frac{\xi}{3},\frac{\delta^*}{2}\}$, $F\in F(G)$ and $\Gamma \subset X_{F,\mathcal{B}(\mu,\xi)}$ be a $(\delta^*,F,\varepsilon^*)$-separated set. Then for any $F^{\prime} \subset F$ with $\frac{|F^{\prime}|}{|F|}> 1- \delta$, $\Gamma$ is a $(\frac{\delta^*}{2},F^{\prime},\varepsilon^*)$-separated set and $\Gamma \subset X_{F^{\prime},B(\mu,2\xi)}$.
\end{lemma}
We remark that the statement of Lemma \ref{eslemma1} is a little bit different from \cite[Lemma 2.6]{zhang2018topological}. But from the proof of \cite[Lemma 2.6]{zhang2018topological}, it is the same thing.

\section{Proofs of main theorems}

The ideas of the proof is to obtain a decomposition for $\bigcup_{n=1}^{\infty}F_n$ which will be used to find suitable orbit segments. Then we construct subsets based on the orbit segments by using the $g$-almost product property and we show those subsets have properties as we want.

First we will recall the quasi-tiling method and a recent result on tilings of amenable groups \cite{downarowicz2019tilings}. The ideas to find the decomposition of $\bigcup_{n=1}^{\infty}F_n$ is from \cite{zhang2018topological}.

Let $\{A_1,\cdots,A_k\} \subset F(G), \alpha >0, \varepsilon>0$ and $A \in F(G)$. The collection $\{A_1,\cdots,A_k\}$ is called {\itshape $\varepsilon$-disjoint} if there exists 
$B_i\subset A_i, i=1,\cdots,k$ with $B_i \cap B_j = \emptyset, 1 \leq i\neq j \leq k$ and $\frac{|B_i|}{|A_i|} > 1-\varepsilon$ for $ i=1,\cdots,k.$ We say $\{A_1,\cdots, A_k\}$ is an $\alpha$-covering of $A$ if $|A\cap \cup_{i=1}^{k} A_i| \geq \alpha|A|$.

We say $\{A_1,\cdots,A_k\} \subset F(G)$ is an {\itshape $\varepsilon$-quasi-tile} of $A$ if there are $C_1,\cdots,C_k\in F(G)$ satisfying:
\begin{enumerate}
	\item $A_iC_i\subset A, i=1,\cdots,k$;
	\item $A_iC_i\cap A_jC_j = \emptyset, 1\leq i \neq j \leq k$;
	\item $\{A_ic \mid c\in C_i\}$ are $\varepsilon$-disjoint for $i=1,\cdots,k$;
	\item $\{A_iC_i\mid i=1,\cdots,k\}$ is a $(1-\varepsilon)$-covering of $A$.
\end{enumerate}  
Such $\{C_1,\cdots,C_k\}$ are called {\itshape tiling centers.}

The following lemma is a fundamental tool in the quasi-tiling theory for amenable groups.

\begin{lemma}[\cite{coornaert2015topological}]\label{tiling1}
	Let $G$ be a countable amenable discrete group and $\{e_G\} \subset F_1\subset F_2\subset\cdots$ be a F{\o}lner sequence. Then for any $\varepsilon \in (0,\frac{1}{4})$ and $N\in \mathbb{N}$, there exist $n_1,n_2,\cdots,n_k\in \mathbb{N}$ and $\delta>0$ with $N\leq n_1<\cdots<n_k$ such that for $F\in F(G)$ which is $(\delta, F_{n_k}F_{n_k}^{-1})$-invariant and $\frac{|F_{n_k}|}{|F|}<\delta,$ then $F$ can be $\varepsilon$-quasi tiled by $F_{n_1},\cdots, F_{n_k}$.
\end{lemma}

\begin{remark}\label{remark1}
	Let $K_1, K_2, \cdots, K_s$ be an $\varepsilon$-quasi-tile of $D \subset G$ and $\{C_1,\cdots,C_s\}$ be the tiling centers. We can modify the tile to get a disjoint $(1-\varepsilon)^2$-covering of $D$ by shrinking every translation of $K_i, \ i = 1, 2, \cdots, s.$ In fact, for each $j$, since $\{K_j c_j \mid c_j\in C_j\}$ are $\varepsilon$-disjoint, we can choose $K_j(c_j) \subset K_j$ with $\frac{|K_j(c_j)|}{|K_j|}\geq 1-\varepsilon$ and the elements in $\{K_j(c_j)c_j\}$ are pairwise disjoint. Thus the elements in the collection $\{K_j(c_j)c_j\mid c_j\in C_j, \ j=1,2,\cdots, s\}$ are pairwise disjoint and 
	\[
	\frac{|\bigcup_{j=1}^{s}\bigcup_{c\in C_j}K_j(c_j)c_j|}{|D|} \geq (1-\varepsilon) \frac{|\bigcup_{j=1}^{s}\bigcup_{c_j\in C_j}K_jc|}{|D|} \geq (1-\varepsilon)^{2}.
	\]
\end{remark}

%{\bf Remark 2.2:} Let $F,K$ be  finite subsets of $G$ and $\gamma >0.$ Suppose $K$ is $(F,\gamma)$-invariant
Also we need a recent tiling result in \cite{downarowicz2019tilings}.

\begin{definition}
	$\mathcal{T}$ is called a {\itshape tiling} of $G$ if there exist a {\itshape shape} set $\mathcal{S} =\{S_j \in F(G) \mid 1\leq j \leq k\}$ and tiling centers $\{C_1, C_2,\cdots, C_k\}$ such that 
	\[
	\mathcal{T}:=\{S_jg \mid g\in C_j, j=1,2,\cdots,k\}
	\]
	with $G = \cup\mathcal{T}$ and $A\cap B =\emptyset$ for $A\neq B \in\mathcal{T}.$ Let $\{\mathcal{T}_k\}_{k\geq 1}$ be a sequence of tilings of $G$, we say $\{\mathcal{T}_k\}_{k\geq 1}$ is {\itshape congruent} if for each $k\geq 1$, each element in $\mathcal{T}_{k+1}$ is a union of elements in $\mathcal{T}_k.$
\end{definition}

The following lemma is part of \cite[Lemma 5.1]{downarowicz2019tilings}.

\begin{lemma}\label{newtiling}
	Fix a positive converging to zero sequence $\{\varepsilon_k\}$ and a sequence $\{K_k\}$ of finite subsets of $G$. There exists a congruent sequence of tilings $\{\tilde{\mathcal{T}_k}\}$ of $G$ such that
	shapes of $\tilde{\mathcal{T}_k}$ are $(K_k,\varepsilon_k)$-invariant.
\end{lemma}

Now we will construct the decomposition of $\bigcup_{n=1}^{\infty}F_n$. Let $\{\beta_k\}_{k=1}^{\infty}$ be a strictly decreasing to 0 positive sequence. We will choose an increasing sequence $M(k)\in \mathbb{N}$ and a sequence of subsets $H_{0}\subset H_1 \subset \cdots \subset G$ in the following way.

\begin{enumerate}
	\item Let $H(0) = \emptyset$. Choose $M(0)\in \mathbb{N}$ such that $F_n$ is $\big(\cup\mathcal{S}_1, \frac{\beta_1}{|\cup\mathcal{S}_1|}\big)$-invariant for any $n\geq M(0)$.
	
	\item Choose $M(1) > M(0)$ such that for any $n \geq M(1)$, $F_n$ is $\big(\cup\mathcal{S}_2, \frac{\beta_2}{|\cup\mathcal{S}_2|}\big)$-invariant. Let $\tilde{F}_1 = \bigcup_{i = M(0) +1}^{M(1)} F_i$, $\tilde{\mathcal{T}_2} =\{ T\in \mathcal{T}_2\mid T\cap \tilde{F}_1 \neq \emptyset\}$
	and $H(1) = \cup \tilde{\mathcal{T}}_2$.
	
	\item Choose $M(2) > M(1)$ such that for any $n \geq M(2)$, $F_n$ is $\big(\cup\mathcal{S}_3, \frac{\beta_3}{|\cup\mathcal{S}_3|}\big)$-invariant and $|H(1)| < \beta_3 |F_n|$. Let $\tilde{F}_2 = \cup_{i = M(0) +1}^{M(2)} F_i$, $\tilde{\mathcal{T}_3} =\{ T\in \mathcal{T}_3\mid T\cap \tilde{F}_2 \neq \emptyset\}$
	and $H(2) = \cup \tilde{\mathcal{T}_3}$.
	
	\item Assume that $M(0) < M(1) < \cdots < M(k-1)$ and $H(0) \subset H(1)\subset \cdots \subset H(k-1)$ have been chosen, then choose $M(k) > M(k-1)$ such that for any $n \geq M(k)$, $F_n$ is 
	$\big(\cup\mathcal{S}_{k+1}, \frac{\beta_{k+1}}{|\cup\mathcal{S}_{k+1|}|}\big)$-invariant and 
	$|H(k-1)| < \beta_{k+1} |F_n|$. Let $\tilde{F}_k = \bigcup_{i = M(0) +1}^{M(k)} F_i$, 
	$\tilde{\mathcal{T}}_{k+1} =\{ T\in \mathcal{T}_{k+1}\mid T\cap \tilde{F}_{k} \neq \emptyset\}$
	and $H(k) = \cup\tilde{\mathcal{T}}_{k+1}$.
	
\end{enumerate}
For $k \geq 1,$ denote
\[
\mathcal{H}_k := \{T\in\mathcal{T}_k \mid T\subset H(k)\setminus H(k-1)\} \text{ and } H_k := \cup \mathcal{H}_k.
\]
Then $H(k)\setminus H(k-1) =H_{k}$ for $k\geq 1.$ 

We call elements in each $\mathcal{T}_k$ {\itshape the standard bricks}. Next we will use standard bricks to cover each $F_n.$ The following lemma shows that most part of $F_n$ can be covered by standard bricks.
\begin{lemma}\cite[Lemma 3.5]{zhang2018topological}\label{standardes1}
	For any $k$ and $M(k-1) < n\leq M(k),$ let 
	\[
	\Lambda_{n}^{1} = \big\{T\in \mathcal{H}_k\mid  T \subset F_n \} 
	\]
	and
	\[
	\Lambda_{n}^{2} = \{T\in \mathcal{H}_{k-1}\mid T\subset F_n\}.
	\]
	Let $\Lambda_{n} = \Lambda_{n}^{1} \cup \Lambda_{n}^{2}$ and $F_n^{\prime} = \cup \Lambda_{n}.$ Then
	\[
	F_{n}^{\prime}\subset F_n \text{ and } |F_{n}^{\prime}| > (1-2\beta_k)|F_n|.
	\]
\end{lemma}

%We mention that we will choose different $\{\mathcal{T}_k\}$ for our use in the following proof.

\subsection{Upper capacity entropy formula}
In this part, we will give the proof of Theorem \ref{mainthm1}.

\begin{proof}
Since $Y_{K} \subset X$, we have $h_{top}^{UC}(Y_K,\{F_n\}) \leq h_{top}(X,\{F_n\}) = h_{top}(X,G)$. Next We will construct a closed subset of $Y_K$ whose upper capacity topological entropy  is close to $h_{top}(X,G)$.

For each $\varepsilon>0$, there exist a finite sequence $\alpha_1,\cdots,\alpha_n$ in $K$ such that each point in $K$ is within $\varepsilon$ of some $\alpha_i.$ As $K$ is connected, by repeating some $\alpha_i$, we can choose this sequence $\alpha_1,\cdots,\alpha_{n}$ so that each point in $K$ is within $\varepsilon$ of some $\alpha_i$ and $D(\alpha_j,\alpha_{j+1}) < \varepsilon$ for each $j.$ Extending this argument, we deduce that there exists a sequence $\{\alpha_j : j =1,2,\cdots\}$ in $K$ so that the closure of
$\{\alpha_j : j>n\}$ for each $n$ equals $K$ and 
\[
\lim\limits_{j\rightarrow\infty}D(\alpha_j,\alpha_{j+1}) = 0.
\]
We define the stretched sequence $\{\alpha^\prime_n\}$ by
\[
\alpha_n^\prime = \alpha_k \text{ if }M(k-1)<n\leq M(k).
\]
The sequences $\{\alpha_n\}$ and $\{\alpha_n^\prime\}$ have the same limit-point set.

	By the variational principle, for any $\eta>0$, there exists $\mu\in E(X,G)$ such that $h_{\mu}(X,G)>h_{top}(X,G) - \eta$. 
	
	Let $\{e_{G}\} \subset K_1\subset K_2 \subset \cdots  $ be a tempered F{\o}lner sequence. Let $\{\xi_k\}_{k=0}^{\infty}$ be a strictly decreasing positive sequence with $\xi_k\rightarrow 0$. By Proposition \ref{sepes1}, there exists $\delta^*>0$ and $\varepsilon^*>0$ such that for the neighborhood $B(\mu,\xi_0)$, there exists $n_{B(\mu,\xi_0),\mu,\eta}^{*}\in \mathbb{N}$ such that  for any $n\geq n_{B(\mu,\xi_0),\mu,\eta}^{*}$, there is a $(\delta^*,K_n,3\varepsilon^*)$-separated set $\Gamma_n \subset X_{K_n,B(\mu,\xi_0)}$ with
\begin{align}
\# \Gamma_n \geq e^{(h_{\mu}(X,G) - \eta)|K_n|}. \label{separatecounting}
\end{align}

%Let $\varepsilon>0$. By Lemma \ref{tiling1}, there exists $n^*\in\mathbb{N}$ such that any $\mathcal{N}\geq n^*$ can be quasi-tiled by some $K_{n_1},\cdots,K_{n_t}$ satisfying $n_{C,\mu,\eta}^{*} \leq n_1 <\cdots<n_t$ with tiling centers $C_{1},\cdots,C_k$. By Remark 2.1, for each $K_{n_j}c$ there exists $K_{n_j}(c) \subset K_{n_j}$ such that $K_{n_j}(c)c$ are pairwise disjoint and $\frac{|K_{n_j}(c)|}{|K_{n_j}|}> 1-\varepsilon$. We mention that we can pick $n^*$ large such that each $K_{N_j}(c)$ is $m(\frac{\varepsilon^*}{4})$-invariant. Let $\Gamma_{n_j}$ be a $(\delta^*,K_{n_j},\varepsilon^*)$-separated set with $\#\Gamma_{n_j} \geq e^{(h_{\mu}(X,G)-\eta)}|K_{n_j|}$. Since $\frac{|K_{n_j}(c)|}{|K_{n_j}|}> 1-\varepsilon$, by Lemma \ref{eslemma1}, we know $\Gamma_{n_j}$ is also $(\frac{\delta^*}{2},K_{n_j}(c),\varepsilon^*)$-separated. For convenience, we will write $\Gamma_{n_j}$ as $\Gamma_{n_j}(c)$ for each $K_{n_j}(c)c$. Define 
%\[
%\Gamma_{F_n} = \prod_{j=1}^{t}\prod_{c\in C_j}\Gamma_{n_j}(c) = \{\vec{x}=(x_{i,j}) \mid \}
%\]
%For each $\vec{x}\in \Gamma_{F_n}$, define
%\[
%B(\delta^*; F_n,\vec(x),\varepsilon^*/4) :=
%\]
%By $g$-almost product property, the set $B(\delta^*; F_n,\vec(x),\varepsilon^*/4)$ is non-empty.

%Thus we have a $(F_n,\frac{\varepsilon^*}{2})$-separated set with cardinality at least $e^{(1-\varepsilon)^2(h_{\mu}(X,G) - \eta)|F_n|}$.

	Let $\{\gamma_k\}_{k=1}^{\infty}$ be a strictly decreasing positive sequence with $\gamma_1 < \min\{\frac{1}{2},\frac{\xi_0}{3},\frac{\delta^*}{2}\}$ and $\gamma_k \rightarrow 0$. By Lemma \ref{tiling1} and Lemma \ref{newtiling}, there exists a sequence of congruent tilings $\mathcal{T}_k$ with shape set $\mathcal{S}_k$ and $N_{1}\leq n_{1,1} <\cdots <n_{1,t_1}<N_{2}\leq n_{2,1}<\cdots<n_{2,t_2}<N_3<\cdots$ such that the following hold:
	Any $S\in \mathcal{S}_k$ can be $\gamma_k$- quasi tiled by $\{K_{n_{k,1}},\cdots,K_{n_{k,t_k}}\}$ with tiling centers $\{C_{k,S,1},\cdots,C_{k,S,t_{k}}\}$ and we denote the $\gamma_k$-quasi-tiling by 
		\begin{align}
		\mathcal{T}_{k,S} = \{K_{c_{k,S,j}}c_{k,S,j}\mid c_{k,S,j}\in C_{k,S,j}, 1\leq j \leq t_k \}. \label{quasitile1}
	\end{align}
Here we assume $N_1 \geq n_{B(\mu,\xi_0),\mu,\eta}^{*}$. By the pointwise ergodic theorem, we can assume that for all $k\geq 1$ and $n\geq N_k$, the set $X_{K_n, B(\alpha_k,\xi_k)}$ is not empty. 
 %By checking the proof of \cite[Theorem 1.1]{ren2020topological}, uniform separation property is only used to calculate entropy. Thus the hypotheses in Theorem \ref{mainthm1} can guarantee the set $Y_K \neq \emptyset$.

From Remark \ref{remark1}, we can modify $\mathcal{T}_{k,S}$ to get a new $\gamma_k$-quasi tile  
\begin{align}
\tilde{\mathcal{T}}_{k,S} = \{\tilde{K}_{c_{k,S,j}}c_{k,S,j}\mid \tilde{K}_{c_{k,S,j}}c_{k,S,j}\subset K_{c_{k,S,j}}c_{k,S,j} \in \mathcal{T}_{k,S}\} \label{qusitile2}\end{align}
such that the following hold:
	\begin{enumerate}
		\item Elements in $\tilde{\mathcal{T}}_{k,S}$ are pairwise disjoint and $\frac{|\tilde{K}_{c_{k,S,j}}|}{|K_{k,j}|}> 1- \gamma_k$;
		\item $\bigcup\tilde{\mathcal{T}}_{k,S} \subset S$ and $|\bigcup \tilde{\mathcal{T}}_{k,S}|>(1-\gamma_k)^2|S|$.
\end{enumerate}
Note that the sequence  obtained from  $\tilde{\mathcal{T}}_{k,S}$
\[
\bigcup_{k\geq 1}\{\tilde{K}_{c_{k,S,j}}\}
\]
is still a F{\o}lner sequence. Let $g:(0,1) \rightarrow (0,1)$ be the mistake-density function as described in the $g$-almost product property and
\[
m:(0,1) \rightarrow F(G) \times (0,1)
\]
be the map(with respect to the mistake-density function $g$) in  Definition \ref{almostspecification}. Choose a decreasing to $0$ sequence $\{\varepsilon_n\}_{n=1}^{\infty}$ with 
\[
\varepsilon_1< \frac{\varepsilon^*}{4} \text{ and } g(\varepsilon_1) < \frac{\delta^*}{4}.
\]
By taking a subsequence of $\bigcup_{k\geq 1}\{\tilde{K}_{c_{k,S,j}}\}$, we can assume for each $n\in\mathbb{N}$, any 
\[
A \in \bigcup_{k\geq n}\left\{ \tilde{K}_{c_{k,S,j}} \mid c_{k,S,j}\in C_{k,S,j}\ 1\leq j\leq t_k, S \in \mathcal{S}_k\right\}
\]
is $m(\varepsilon_n)$-variant.

Let  $\mathcal{N}$ be any integer larger than $M(1)$. Then there exists $\kappa\geq 2$ such that $M(\kappa -1) < \mathcal{N} \leq M(\kappa)$. 

For $k=\kappa-1, \kappa$ and $S\in \mathcal{S}_k$, let $\mathcal{T}_{k,S}$ and $\tilde{\mathcal{T}}_{k,S}$ be as in (\ref{quasitile1}) and (\ref{qusitile2}), respectively. For $K_{c_{k,S,j}}c_{k,S,j} \in \mathcal{T}_{k,S}$, there is a $(\delta^*,K_{c_{k,S,j}},3\varepsilon^*)$-separated set $\Gamma_{c_{k,S,j}} \subset X_{K_{c_{k,S,j}},B(\mu,\xi_0)}$ with $\# \Gamma_{c_{k,S,j}} \geq e^{|K_{c_{k,S,j}}|(h_{\mu}(X,G) - \eta)}.$ By Lemma \ref{eslemma1}, the set $\Gamma_{c_{k,S,j}}$ is also $(\delta^*/2,\tilde{K}_{c_{k,S,j}},3\varepsilon^*)$-separated and contained in $X_{\tilde{K}_{c_{k,S,j}},B(\mu,2\xi_0)}.$

Define 
\begin{align*}
\Gamma(S) & := \prod_{\tilde{K}_{c_{k,S,j}}c_{k,S,j}\in \tilde{\mathcal{T}}_{k,S}}\Gamma_{{c_{k,S,j}}}\\
	& := \{\vec{x} = (x_{c_{k,S,j}}) \mid x_{c_{k,S,j}} \in \Gamma_{{c_{k,S,j}}}\}.
\end{align*}
For each $\vec{x}\in \Gamma(S)$, denote
\begin{align}
B(g;S,\vec{x},\varepsilon_k) := \bigcap_{\tilde{K}_{k,S,j}c_{k,S,j}\in \tilde{\mathcal{T}}_{k,S}}B(g;\tilde{K}_{c_{k,S,j}}c_{k,S,j},c_{k,S,j}^{-1}x_{c_{k,S,j}},\varepsilon_k), \label{sball}
\end{align}
where the dynamical ball $B(g;\tilde{K}_{c_{k,S,j}}c_{k,S,j},c_{k,S,j}^{-1}x_{c_{k,S,j}},\varepsilon_k)$ is defined in Definition \ref{almostspecification}. By the $g$-almost product property, $B(g;S,\vec{x},\varepsilon_k) \ne \emptyset.$

For $k\geq \kappa+1$ and $S\in\mathcal{S}_k$, since the set $X_{K_n,B(\alpha_k,\gamma_k)}$ is not empty for $n\geq N_k$, we pick $x_{c_{k,S,j}} \in X_{K_{c_{k,S,j}}, B(\alpha_k,\xi_k)}$ for $K_{c_{k,S,j}}c_{k,S,j} \in \mathcal{T}_{k,S}$. By a simple calculation, we have  $x_{c_{k,S,j}} \in X_{\tilde{K}_{c_{k,S,j}}, B(\alpha_k,\xi_k+2\gamma_k)}.$ Let $\Gamma_{c_{k,S,j}} = \{x_{c_{k,S,j}}\}$.

Let

	\begin{align*}
	\Gamma(S) := \prod_{\tilde{K}_{c_{k,S,j}}c_{k,S,j}\in \tilde{T}_{k,S}}x_{c_{k,S,j}} = \{\vec{x}=(x_{c_{k,S,j}})\mid x_{c_{k,s,j}} \in \Gamma_{c_{k,S,j}}, \tilde{K}_{c_{k,S,j}}c_{k,S,j}\in \tilde{T}_{k,S}\}.
	\end{align*}
For each $\vec{x}\in \Gamma(S)$, we denote
\[
B(g;S,\vec{x},\varepsilon_k) := \bigcap_{\tilde{K}_{c_{k,S,j}}c_{k,S,j}\in \tilde{\mathcal{T}}_{k,S}}B(g;\tilde{K}_{c_{k,S,j}}c_{k,S,j},c_{k,S,j}^{-1}x_{c_{k,S,j}},\varepsilon_k).
\]
By the $g$-almost product property, $B(g;S,\vec{x},\varepsilon_k) \ne \emptyset.$

For any integer $q >0$, let 
\[
Y_{q}^{(\mathcal{N})} := \bigcap_{i=\kappa-1}^{\kappa+q}\bigcup_{Sd \in H_{i}}\bigcup_{\vec{x}\in \Gamma(S)} B(g;Sd,d^{-1}\vec{x},\varepsilon_i),
\]
where $B(g;Sd,d^{-1}\vec{x},\varepsilon_i) := \bigcap_{\tilde{K}_{k,S,j}c_{k,S,j}\in \tilde{\mathcal{T}}_{k,S}}B(g;\tilde{K}_{c_{k,S,j}}c_{k,S,j}d,d^{-1}c_{k,S,j}^{-1}x_{c_{k,S,j}},\varepsilon_i). $
By the $g$-almost product property, the set $Y_{q}^{(\mathcal{N})}$ is not empty. Let $Y^{(\mathcal{N})} = \bigcap_{q\geq 1}Y_{q}^{(\mathcal{N})}$,
then $Y^{(\mathcal{N})}$ is not empty.

Next we will show some properties of  $Y^{(\mathcal{N})}$:

\begin{enumerate}
	\item $Y^{(\mathcal{N})}\subset Y_{K}$;
	\item There exists $Y\subset Y^{(\mathcal{N})}$ which is $(F_{\mathcal{N}},\varepsilon^*)$-separated and\\ $\# Y \geq e^{(1-3\eta)|F_{\mathcal{N}}|(h_{top}(X,G)-2\eta)}$.
\end{enumerate}

{\bf Proof of Item (1):} 
Choose $y\in Y^{(\mathcal{N})}$ and $F_n.$ Suppose $n$ large such that $M(k-1)<n\leq M(k)$ and $k \geq \kappa+2.$ Let $\Lambda_{n}^{(1)}$ and $\Lambda_{n}^{(2)}$ be as described in Lemma \ref{standardes1}. Using the triangle inequality, Lemma \ref{eslma1} and Lemma \ref{eslemma1},

\begin{align}
D(\mathcal{E}_{F_n}(y),\alpha_n^{\prime}) &\leq D(\mathcal{E}_{F_n^{\prime}}(y),\alpha_k) + D(\mathcal{E}_{F_n}(y),\mathcal{E}_{F_n^{\prime}}(y)) \notag\\
&\leq 4\beta_k + \sum_{Sd\in \Lambda_{n}^1}\frac{|Sd|}{|F_n^{\prime}|}D(\mathcal{E}_{Sd}(y),\alpha_k) + \sum_{Sd\in \Lambda_{n}^2}\frac{|Sd|}{|F_n^{\prime}|}D(\mathcal{E}_{Sd}(y),\alpha_k)  \\
&\leq 4\beta_k + \sum_{Sd\in \Lambda_{n}^1}\frac{|Sd|}{|F_n^{\prime}|}\left(D(\mathcal{E}_{Sd}(y),\mathcal{E}_{\tilde{S}d}(y))+D(\mathcal{E}_{\tilde{S}d}(y),\alpha_k)\right) \notag \\
&+ \sum_{Sd\in \Lambda_{n}^2}\frac{|Sd|}{|F_n^{\prime}|}\left(D(\mathcal{E}_{Sd}(y),\mathcal{E}_{\tilde{S}d}(y))+D(\mathcal{E}_{\tilde{S}d}(y),\alpha_k)\right)\\
&\leq 4\beta_k + \sum_{Sd\in \Lambda_{n}^1}\frac{|Sd|}{|F_n^{\prime}|}\left(4\gamma_k+D(\mathcal{E}_{\tilde{S}d}(y),\alpha_k)\right)\notag\\
&+ \sum_{Sd\in \Lambda_{n}^2}\frac{|Sd|}{|F_n^{\prime}|}\left(4\gamma_{k-1}+D(\mathcal{E}_{\tilde{S}d}(y),\alpha_k)\right)\\
&\leq 4\gamma_{k-1} + 4\beta_k + \sum_{Sd\in \Lambda_{n}^1}\frac{|Sd|}{|F_n^{\prime}|} (\xi_k + \varepsilon_k + g(\varepsilon_k))\notag \\
&+ \sum_{Sd\in \Lambda_{n}^2}\frac{|Sd|}{|F_n^{\prime}|}(\xi_{k-1} + \varepsilon_{k-1} + g(\varepsilon_{k-1}) + D(\alpha_{k-1},\alpha_k))\\
&\leq 4\gamma_{k-1} + 4\beta_k + \xi_{k-1} + \varepsilon_{k-1} + g(\varepsilon_{k-1}) + D(\alpha_{k-1},\alpha_k). \label{eslemma3}
\end{align}

This implies $\lim\limits_{n\rightarrow \infty}D(\mathcal{E}_{F_n}(y),\alpha_n^\prime) = 0$. Hence we have $Y^{(\mathcal{N})}\subset Y_{K}.$

{\bf Proof of Item (2):} Define 
\begin{align*}
\Gamma(F_{\mathcal{N}}) :=& \prod_{Sd\in \Lambda_{\mathcal{N}}}\Gamma(S)\\
=&\{\vec{x} = (\vec{x}_{S})\mid \vec{x}_{S}\in \Gamma(S), Sd\in \Lambda_{\mathcal{N}}\}\\
=& \{\vec{x} = (x_{c_{k,S,j}})  \mid x_{c_{k,S,j}}\in \Gamma_{c_{k,S,j}},\  K_{c_{k,S,j}}c_{k,S,j}\in \mathcal{T}_{k,S},\  Sd\in \Lambda_{\mathcal{N}} \}.
\end{align*}

From the definition, we have
\begin{align}
\# \Gamma(F_{\mathcal{N}}) = \prod_{Sd\in \Lambda_{\mathcal{N}}} \# \Gamma(S). \label{caless}
\end{align}

Write $\Gamma(F_{\mathcal{N}}) = \{\vec{x}^1,\cdots, \vec{x}^r \}$, where $r=\#\Gamma(F_{\mathcal{N}}).$ For $\vec{x}^u \neq \vec{x}^v \in \Gamma(F_{\mathcal{N}})$, there exists $Sd\in  \Lambda_{\mathcal{N}}$ such that $\vec{x}^{u}_{S} \neq \vec{x}^{v}_{S}$ which implies there exists $\tilde{K}_{c_{k,S,j}}c_{k,S,j} \in \tilde{\mathcal{T}}_{k,S}$ such that $x^{u}_{c_{k,S,j}}$ and $x^{v}_{c_{k,S,j}}$ are $(\delta^*,\tilde{K}_{c_{k,S,j}},3\varepsilon^*)$-separated.

For each $\vec{x}^{u} \in \Gamma(F_{\mathcal{N}})$, there is $x^{u}\in Y^{(\mathcal{N})}$ such that $x^{u}\in B(g;Sd, d^{-1}\vec{x}_{S}, \tilde{\varepsilon})$ for every $Sd\in \Lambda_{\mathcal{N}}$, where $\tilde{\varepsilon} = \varepsilon_\kappa$ if $Sd\in \Lambda_{\mathcal{N}}^1$ and  $\tilde{\varepsilon} = \varepsilon_{\kappa-1}$ if $Sd\in \Lambda_{\mathcal{N}}^2$. Denote $Y=\{x^{u}\mid \vec{x}^{u}\in \Gamma(F_{\mathcal{N}})\}$. Thus there is a one-to-one map $\Phi$ from $\Gamma(F_{\mathcal{N}})$ to $Y$.

To this end, fix $x^{u} \neq x^{v} \in Y$. Let $\vec{x}^{u}=(x^u_{c_{k,S,j}}) = \Phi^{-1}(x^{u}) \in \lambda_{\mathcal{N}} $ and $\vec{x}^{v}=(x^{v}_{c_{k,S,j}})=\Phi^{-1}(x^{v}) \in \Gamma(F_{\mathcal{N}})$. Then there is $Sd\in\Lambda_{\mathcal{N}}$ and $\tilde{K}_{c_{k,S,j}}c_{k,S,j} \in \tilde{\mathcal{T}}_{k,S}$ such that  $x^{u}_{c_{k,S,j}}$ and $x^{v}_{c_{k,S,j}}$ are $(\delta^*,\tilde{K}_{c_{k,S,j}},3\varepsilon^*)$-separated. And also $x^{u}\in B(g;Sd, d^{-1}\vec{x}^{u}_{S}, \tilde{\varepsilon}), x^{v}\in B(g;Sd, d^{-1}\vec{x}^{v}_{S}, \tilde{\varepsilon})$  implies 
\[
\frac{|\{s\in \tilde{K}_{c_{k,S,j}}\mid \rho(sx^u,sx^{u}_{c_{k,S,j}}) > {\tilde\varepsilon}\}|}{|\tilde{K}_{c_{k,S,j}}|} > \delta^* \text{ and } \frac{|\{s\in \tilde{K}_{c_{k,S,j}}\mid \rho(sx^v,sx^{v}_{c_{k,S,j}}) > {\tilde\varepsilon}\}|}{|\tilde{K}_{c_{k,S,j}}|} > \delta^*.
\]
Then there is $s\in \tilde{K}_{c_{k,S,j}}$ such that 
\[
\rho(sx^u,sx^v) \geq 3\varepsilon^* - 2\tilde{\varepsilon} > \varepsilon^*.
\]
Based on items (1) and (2), we get that $Y^{(\mathcal{N})}\subset Y_{K}$ and for each $\eta>0$, $h_{top}^{UC}(Y_K,\{F_n\}) \geq h_{top}^{UC}(Y^{\mathcal{(N)}},\{F_n\}) > h_{\mu}(X,G) >h_{top}(X,G) - 2\eta.$ By the arbitrariness of $\eta$, we have 
$h_{top}^{UC}(X,G) \geq h_{top}(X,G)$, which finish the proof of Theorem \ref{mainthm1}.

\end{proof}

\subsection{Packing entropy formula}

In this part, we will proceed the proof of Theorem \ref{mainthm2}.

First, we will construct a closed subset of $Y_K$ whose packing topological entropy is close to $\sup\{h_{\mu}(X,G)\mid \mu\in K\}$. Let $\eta>0$ and $H^{*}=\sup\{h_{\mu}(X,G)\mid \mu\in K\} - \eta$ and $h^* = \inf\{h_{\mu}(X,G)\mid \mu\in K\}-\eta.$ Pick  $\mu_{\max}\in K$ such that $h_{\mu_{max}}(X,G)> H^*.$ Since $K$ is convex and compact, for each $n\in \mathbb{N},$ we can choose $\{\mu_{n,1},\cdots,\mu_{n,\mathcal{M}_n}\} \subset K$ such that
\[
K\subset \bigcup_{i}\mathcal{B}(\mu_{n,i},\frac{1}{n}), D(\mu_{n,i},\mu_{n,i+1}) \leq \frac{1}{n} \text{ and } \mu_{n,\mathcal{M}_{n-1}} = \mu_{n,\mathcal{M}_n} = \mu_{\max}.
\]
Here we repeat $\mu_{max}$ twice since two levels of standard bricks are needed to cover each $F_n$. Denote $\{\alpha_j\}=\{\mu_{1,1},\mu_{1,2},\cdots,\mu_{1,\mathcal{M}_1},\mu_{2,1},\cdots,\mu_{2,\mathcal{M}_2},\cdots\}$.
Then $K = \overline{\{\alpha_j\}}_{j\geq n}$ for each $n\geq 1$ and $D(\alpha_j,\alpha_{j+1}) \rightarrow 0$ as $j\rightarrow \infty.$ Define the stretched sequence $\{\alpha^\prime_n\}$ by
\[
\alpha_n^\prime = \alpha_k \text{ if }M(k-1)<n\leq M(k).
\]
The sequences $\{\alpha_n\}$ and $\{\alpha_n^\prime\}$ have the same limit-point set. 

Let $\{\xi_k\},\{\eta_k\}$ and $\{\gamma_k\}$ be sequences of  strictly decreasing positive numbers with $\gamma_1 < \min\{\frac{1}{2},\frac{\xi_0}{3},\frac{\delta^*}{2}\},\eta_1 <\eta$ and $\xi_k\rightarrow 0, \eta_k \rightarrow0, \gamma_k \rightarrow 0$. Let $\{K_n\}$ be a tempered F{\o}lner sequence with $e_{G}\in K_1 \subset K_2 \subset \cdots.$ For the given $\eta>0$, by Lemma \ref{separatinges1}, we can find $\delta^*>0$ and $\varepsilon^*>0$ such that for $\mu\in M(X,G)$ and any neighborhood $C\subset M(X)$ of $\mu$, there exists $n^{*}_{C,\mu,\eta}$ such that for $n\geq n^{*}_{C,\mu,\eta}$,
\begin{align}
N(C;\delta^*,K_n,\varepsilon^*) \geq e^{|K_n|(h_{\mu}(X,G)-\eta)}. \label{number2}
\end{align}

 By Lemma \ref{tiling1} and Lemma \ref{newtiling}, there exists a sequence of congruent tilings $\mathcal{T}_k$ with shape set $\mathcal{S}_k$ and $N_{1}\leq n_{1,1} <\cdots <n_{1,t_1}<N_{2}\leq n_{2,1}<\cdots<n_{2,t_2}<N_3<\cdots$ such that the following hold:
Any $S\in \mathcal{S}_k$ can be $\gamma_k$- quasi tiled by $\{K_{k,1},\cdots,K_{k,t_k}\}$ with tiling centers $\{C_{k,S,1},\cdots,C_{k,S,t_{k}}\}$ and we denote the $\gamma_k$-quasi-tiling by 
\begin{align}
\mathcal{T}_{k,S} = \{K_{{k,j}}c_{k,S,j}\mid c_{k,S,j}\in C_{k,S,j}, 1\leq j \leq t_k \}. \label{quasitiling3}
\end{align}
We also assume that  $N_k \geq n^*_{B(\alpha_k,\xi_k),\alpha_k,\eta_k}$ for each $k$.

From Remark \ref{remark1}, we can modify $\mathcal{T}_{k,S}$ to get a new $\gamma_k$-quasi tile  
\begin{align}
\tilde{\mathcal{T}}_{k,S} = \{\tilde{K}_{c_{k,S,j}}c_{k,S,j}\mid K_{c_{k,S,j}}c_{k,S,j} \in \mathcal{T}_{k,S}\} \label{qusitile4}\end{align}
such that the following hold:
\begin{enumerate}
	\item Elements in $\tilde{\mathcal{T}}_{k,S}$ are pairwise disjoint and $\frac{|\tilde{K}_{c_{k,S,j}}|}{|K_{c_{k,S,j}}|}> 1- \gamma_k$;
	\item $\bigcup\tilde{\mathcal{T}}_{k,S} \subset S$ and $|\bigcup \tilde{\mathcal{T}}_{k,S}|>(1-\gamma_k)^2|S|$.
\end{enumerate}
From (\ref{number2}), we get the existence of a $(\delta^*,K_{c_{k,S,j}},\varepsilon^*)$-separated set $\Gamma_{c_{k,S,j}} \subset X_{K_{c_{k,S,j}},B(\alpha_k,\xi_k)}$ with $\#\Gamma_{c_{k,S,j}} \geq e^{|K_{c_{k,S,j}}|(h_{\alpha_k}(X,G) - \eta_k)}$.   By Lemma \ref{eslemma1}, the set $\Gamma_{c_{k,S,j}}$ is also $(\delta^*/2,K_{c_{k,S,j}},\varepsilon^*)$-separated and contained in $X_{\tilde{K}_{c_{k,S,j}},B(\alpha_k,2\xi_k)}.$

Let
\begin{align*}
\Gamma(S) & := \prod_{\tilde{K}_{c_{k,S,j}}c_{k,S,j}\in \tilde{\mathcal{T}}_{k,S}}\Gamma_{{c_{k,S,j}}}\\
& = \{\vec{x} = (x_{c_{k,S,j}}) \mid x_{c_{k,S,j}} \in \Gamma_{{c_{k,S,j}}}, K_{c_{k,S,j}}c_{k,S,j}\in \mathcal{T}_{k,S}\}.
\end{align*}
For $Sd\in \mathcal{T}_k,$ let $\Gamma(Sd) = \Gamma(S)$. For each $\vec{x}\in \Gamma(S)$, denote
\[
B(g;S,\vec{x},\varepsilon_k) := \bigcap_{\tilde{K}_{c_{k,S,j}}c_{k,S,j}\in \tilde{\mathcal{T}}_{k,S}}B(g;\tilde{K}_{c_{k,S,j}}c_{k,S,j},c_{k,S,j}^{-1}x_{c_{k,S,j}},\varepsilon_k).
\]

For  $q\geq 1$, let 
\[
Y^{q}_{\eta} := \bigcap_{j=0}^{q}\big(\bigcup_{Sd \in H_{j}}\bigcup_{\vec{x}\in \Gamma(S)}B(g;Sd,d^{-1}\vec{x},\varepsilon_k)\big),
\]
where $B(g;Sd,d^{-1}\vec{x},\varepsilon_k) = B(g;S,\vec{x},\varepsilon_k).$
By the $g$-almost product property, the set $Y^{q}_{\eta}$ is a non-empty closed set.

For each $q \geq 1$, define
\begin{align}
\Gamma_{q} &:= \prod_{j=1}^{q}\prod_{Sd\in H_j}\Gamma(Sd)\notag\\
&= \{\vec{x} = (\vec{x}_{j,Sd}) \mid \vec{x}_{j,Sd} \in \Gamma(Sd), Sd\in H_j \text{ and } j=1,2,\cdots,q\}.
\end{align}
For each $\vec{x}_q \in \Gamma_q$, pick $x_q\in X$ such that $x_q\in B(g;Sd, d^{-1}x_{j,Sd})$ for all $Sd\in H_j, j=1,2,\cdots,q.$ Let 
\[
Z_{q} := \{x_q \mid \vec{x}_q \in \prod_{j=1}^{q}\prod_{Sd\in H_j}\Gamma(Sd)\}. 
\]
To this end, define
\[
\mu_{q} := \frac{1}{\#Z_q} \sum_{x\in Z_q}\delta_x,
\]
where $\delta_x$ is the Dirac measure at $x$. For any fixed $l\geq 0$ and $p\geq 0$, we have $\mu_{l+p}(Y^{l}_{\eta}) =1$ since $Y^{l}_{\eta} \supset Y^{l+p}_{\eta}$ and
$\mu_{l+p}(Y^{l+p}_{\eta}) =1$. Let $\mu$ be a limit measure of $\{\mu_q\}$ i.e. there exists a sequence $\{q_n\}$ such that $\mu_{q_n} \rightarrow \mu.$ Since $Y^{l}_{\eta}$ is closed, we have $\mu(Y^{l}_{\eta}) \geq \limsup\limits_{n\rightarrow \infty}\mu_{q_n}(Y^{l}_{\eta}) =1$. It follows that $\mu(Y_{\eta}) = 1.$

Take $x\in Y_{\eta}$. For $i\geq 1$, let $L_{i}=\sum_{j=1}^{i}\mathcal{M}_j$. Then for $q_n \geq L_{i}$, one has
\begin{align*}
\mu_{q_n}(B_{F_{M(L_i)}}(x,\varepsilon)) & \leq \frac{1}{\prod_{Sd\in \Lambda_{M(L_i)}}\# \Gamma(S)}\\
& = \frac{1}{\prod_{Sd\in \Lambda_{M(L_i)}} \prod_{K_{c_{k,S,j}c_{k,S,j} \in \mathcal{T}_{k,S}}}\#\Gamma_{{c_{k,S,j}}}}\\
&\leq \frac{1}{\prod_{Sd\in \Lambda_{M(L_i)}} \prod_{K_{c_{k,S,j}c_{k,S,j} \in \mathcal{T}_{k,S}}}e^{|\tilde{K}_{c_{k,S,j}}|(h_{\mu_{max}(X,G)} - \eta)}}\\
&\leq e^{-(1-\gamma_{L_{i-1}})^2(1-2\beta_{L_{i-1}})|F_{M(L_i)}|(h_{\mu_{max}(X,G)} - \eta)}.
\end{align*} 
Therefore $\mu(B_{F_{M(L_i)}}(x,\varepsilon)) \leq e^{-(1-\gamma_{L_{i-1}})^2(1-2\beta_{L_{i-1}})|F_{M(L_i)}|(h_{\mu_{max}(X,G)} - \eta)}$.
Thus we have 
\begin{align*}
\overline{h}^{loc}_{\mu}(x,\varepsilon,\{F_n\}) &\geq \limsup_{i\rightarrow \infty}-\frac{1}{|F_{M(L_i)}|}\log \mu\left(B_{F_{M(L_i)-1}}(x,\varepsilon)\right)\\
&\geq \limsup\limits_{i\rightarrow \infty}{{(1-\gamma_{L_{i-1}})^2(1-2\beta_{L_{i-1}})(h_{\mu_{max}(X,G)} - \eta)}}\\
&\geq H^* - \eta,
\end{align*}
where $\overline{h}^{loc}_{\mu}(x,\varepsilon,\{F_n\})$ is defined in Definition \ref{localentropy}.

By Theorem \ref{packingvp}, we have $h_{top}^{P}(Y_{K},\{F_n\}) \geq h_{top}^{P}(Y_{\eta},\{F_n\}) \geq H^* -\eta = \sup\{h_{\mu}(X,G) \mid \mu\in K \} - 2\eta$. By the arbitrariness of $\eta$, we have $h_{top}^{P}(Y_{K},\{F_n\}) \geq \sup\{h_{\mu}(X,G) \mid \mu\in K \}.$

Next we will show $h_{top}^{P}(Y_K,\{F_n\}) \leq \sup\{ h_{\mu}(X,G)\mid \mu\in K\}$. Recall that $V(x,\{F_n\})$ is the limit-point set of  $\{\mathcal{E}_{F_n}(x)\}.$ Define
\[
Y^{K} := \{y\in X \mid V(y,\{F_n\}) \cap K \neq \emptyset\}.
\]
From the definition, $Y_{K} \subset Y^{K}$.
For $\delta>0$ and $n\in \mathbb{N}$, set 
\[
\mathcal{R}(K,\delta,n) := \{x\in X \mid \mathcal{E}_{F_n}(x) \in B(K,\delta)\},
\]
where $B(K,\delta) = \{\nu\in M(X) \mid \exists\  \mu\in K\text{ such that } D(\mu,\nu) \leq \delta\}$.

Fix $\varepsilon>0$ and let $N(K,\delta,F_n,\varepsilon)$ denote the smallest number of balls $B_{F_n}(x,\varepsilon)$ required to cover $\mathcal{R}(K,\delta,n)$. Notice that $N(K,\delta,F_n,\varepsilon)$ does not increase as $\delta$ decreases and does not decrease as $\varepsilon$ decreases. As a result, the following limit denoted by $\Theta(Y^{K},\{F_n\})$ exists
\begin{align}
	\Theta(Y^{K},\{F_n\}) := \lim\limits_{\varepsilon \rightarrow 0}\lim\limits_{\delta \rightarrow 0} \limsup_{n\rightarrow \infty}\frac{1}{|F_n|}\log N(K,\delta,F_n,\varepsilon) \label{entropydef1}.
\end{align}

Let $R(K,\delta,k) = \bigcap_{n=k}^{\infty} \mathcal{R}(K,\delta,n) = \bigcap_{n=k}^{\infty}\{y\in X \mid \mathcal{E}_{F_n}(x)\in B(K,\delta)\}$. Let $\tilde{N}(K,\delta,F_n,\varepsilon)$ be the smallest number of balls $B_{F_n}(x,\varepsilon)$ used to cover $R(K,\delta,n)$. From the definition, \begin{align}
\tilde{N}(K,\delta,F_n,\varepsilon) \leq N(K,\delta,F_n,\varepsilon). \label{comparejust1}
\end{align}
For each $\delta>0$, we have
\begin{align}
	Y^{K} &= \{x\in X \mid \lim\limits_{n\rightarrow\infty} D(\mathcal{E}_{F_n}(x),K) = 0\}\notag\\
	&\subset \bigcup_{k=1}^{\infty}R(K,\delta,k).
\end{align}

For each $k$, let $R(K,k) = \bigcap_{\delta>0}R(K,\delta,k).$ Then $Y^{K} \subset \bigcup_{k=1}^{\infty}R(K,k)$. By Proposition \ref{packingproperty}(2), we have 
\begin{align}
h_{top}^{P}(Y^{K},\{F_n\}) \leq \sup_{k\geq 1}h_{top}^{P}(R(K,k)). \label{essect30}
\end{align}
By (\ref{comparejust1}), we have
\begin{align}
h_{top}^{UC}(R(K,k),\{F_n\}) \leq \Theta(Y^{K},\{F_n\}). \label{comparejust2}
\end{align}
By Proposition \ref{packinguc} and (\ref{essect30} - \ref{comparejust2}),  $h_{top}^{P}(Y^{K},\{F_n\}) \leq \Theta(Y^{K},\{F_n\}).$ By Proposition \ref{packingproperty}, we have
\begin{align}
	h_{top}^{P}(Y_{K},\{F_n\}) \leq \Theta(Y^{K},\{F_n\}). \label{sect3part1}
\end{align}

We just need to show that
\[
\Theta(Y^{K},\{F_n\}) \leq \sup_{\mu\in K}h_{\mu}(X,G).
\]

By (\ref{entropydef1}), for any $\gamma>0,$ there exists $\varepsilon_0$ small such that for all $0<\varepsilon<\varepsilon_0$ the following holds
\begin{align}
\lim\limits_{\delta\rightarrow 0}\limsup_{n\rightarrow \infty} \frac{1}{|F_n|} \log N(K,\delta,F_n,\varepsilon) > \Theta(G^K,\{F_n\}) - \frac{\gamma}{3}. \label{essec31}
\end{align}
For each $k\geq 1$, set $\varepsilon_k = \frac{\varepsilon_0}{2^k}$. From (\ref{essec31}), there is $\delta_k$ samll with $\delta_k\rightarrow 0$ and 
\begin{equation}
	\limsup_{n\rightarrow \infty} \frac{1}{|F_n|}\log N(K,\delta_k,F_n,\varepsilon_k) \geq \Theta(Y^{K},\{F_n\}) - \frac{2\gamma}{3}. \label{essect32}
\end{equation} 
Choose $n_k\in\mathbb{N}$ such that
\begin{equation}
	N(K,\delta_k,F_{n_k},\varepsilon_k) \geq e^{|F_{n_k}|(\Theta(Y^{K},\{F_n\}) - \gamma)}. \label{essect33}
\end{equation}
Let $C_k$ be the centers of one covering of $\mathcal{R}(K,\delta_k,n_k)$ with balls $B_{F_{n_k}}(x,\varepsilon_k)$ and $\# C_k = N(K,\delta_k,F_{n_k},\varepsilon_k).$ Let $\mu_n\in M(X)$ be defined by $\mu_{n} = \frac{1}{\# C_k }\sum_{x\in C_k}\mathcal{E}_{F_{n_k}}(x)$. Let $\mu\in M(X)$ be a limit point of the sequence $\{\mu_k\}$. Since $\{F_n\}$ is a F{\o}lner sequence, the measure $\mu$ is in $M(X,G)$. For each $x\in C_k$, we can pick one point $y(x)\in \mathcal{R}(K,\delta_k,n_k)\cap B_{F_{n_k}}(x,\varepsilon_k)$. Then $D(\mathcal{E}_{F_{n_k}}(x),K) \leq D(\mathcal{E}_{F_{n_k}}(y(x)),K) + \delta_k \leq \varepsilon_k + \delta_k$. As a result we have $D(\mu,K) = 0$. Since $K$ is closed, $\mu\in K.$

By Lemma \ref{upperlemma}, 
\begin{equation}
	h_{\mu}(X,G) \geq \limsup_{k\rightarrow \infty} \frac{1}{|F_{n_k}|} \log N(K,\delta_k,F_{n_k},\varepsilon_k). \label{essect34}
\end{equation}

Combining (\ref{essect33}) and (\ref{essect34}), we have $h_{\mu}(X,G) \geq \Theta(Y^{K},\{F_n\}) - \gamma$.  By the arbitrariness of $\gamma$, we have
\begin{align}
	h_{\mu}(X,G) \geq \Theta(Y^{K},\{F_n\}). \label{sect3part2}
\end{align}
We get $h_{top}^{P}(G_K,\{F_n\}) \leq \sup_{\mu\in K} h_{\mu}(X,G)$ by combining (\ref{sect3part1}) and (\ref{sect3part2}).

\section{An application in multi-fractal analysis}
 We present one application  of the above results.
 
Let $(X,G)$ be a dynamical system and $\{F_n\}$ be a F{\o}lner sequence. Let $\varphi: X\rightarrow \mathbb{R}$ be a continuous function. Denote by 
\[
L_{\varphi} = [\inf_{\mu\in M(X,G)}\int \varphi \dd\mu, \sup_{\mu\in M(X,G)}\int \varphi \dd\mu].
\]
 For any $a\in L_{\varphi}$, define the level set with respect to $\{F_n\}$ 
 \[
 R_{\varphi}(a) = \{x\in X \mid \lim\limits_{n\rightarrow\infty}\frac{1}{|F_n|}\sum_{s\in F_n}\varphi(sx) = a\}.
 \]
 
 For $\mu\in M(X,G)$, let $Y_{\mu} =\{x\in X\mid V(x,\{F_n\}) = \{\mu\}\}$ be the saturated set with respect to $\{\mu\}$.
 
\begin{proposition}
	Let $(X,G)$ be a system satisfying the $g$-almost product property and $\{F_n\}$ be a F{\o}lner sequence %satisfying $\frac{|F_n|}{|\log n|} \rightarrow \infty$. 
	Let $\varphi: X \rightarrow \mathbb{R}$ be a continuous function. Then for any $a\in L_{\varphi}$, we have
	\[
	h_{top}^{UC}(R_{\varphi}(a),\{F_n\}) = h_{top}(X,G).
	\]
	
	Moreover, if the system also has the uniform separation property and the sequence $\{F_n\}$ satisfies $\frac{|F_n|}{|\log n|} \rightarrow \infty$, we have
	\[
	h_{top}^{P}(R_{\varphi}(a),\{F_n\}) = \sup\{h_{\mu}(X,G) \mid \mu\in M(X,G), \int \varphi \dd \mu = a\}.
	\]
\end{proposition}

\begin{proof}
	Let $F(a)=\{\nu\in M(X,G) \mid \int \varphi \dd\nu = a\}$. It is clear that for any $\nu\in F(a)$, $Y_{\nu}\subset R_{\varphi}(a)$. So 
	\[
	h_{top}^{UC}(Y_{\nu},\{F_n\}) \leq h_{top}^{UC}(R_{\varphi}(a),\{F_n\}),
	 \ \forall \nu\in F(a).
	\]
	Combining Theorem \ref{mainthm1}, we have $h_{top}^{UC}(R_{\varphi}(a),\{F_n\}) = h_{top}(X,G)$.
	
	Now we prove the moreover part. It is obvious that 
	\[
	R_{\varphi}(a) = \{x\in X\mid V(x,\{F_n\})  \subset F(a)\}.
	\]
	
	With the similar discussion in the proof of the upper bound of $h_{top}^{P}(Y_K,\{F_n\})$, we get that $h_{top}^{P}(R_{\varphi}(a),\{F_n\}) \leq \sup\{h_{\mu}(X,G) \mid \mu\in F(a)\}$. For any $\mu\in F(a)$, one has $Y_{\mu}\subset R_{\varphi}(a)$ which implies $h_{top}^{P}(Y_{\mu},\{F_n\})\leq h_{top}^{P}(R_{\varphi}(a),\{F_n\}).$ By Theorem \ref{mainthm2}, we have $h_{top}^{P}(Y_{\mu},\{F_n\}) = h_{\mu}(X,G)$. Then we have $h_{top}^{P}(R_{\varphi}(a),\{F_n\}) = \sup\{h_{\mu}(X,G) \mid \mu\in M(X,G), \int \varphi \dd \mu = a\}$
\end{proof}

\section{Appendix}

Before proceeding the proof, we give the definition of $f$-neighborhood.

\begin{definition}
	
	An $f$-neighborhood of $\mu\in M(X)$ is the set of the form 
	\[
	F^{(\alpha)} := \{\nu\in M(X) \mid \big|\langle f_i,\mu\rangle - \langle f_i,\nu\rangle\big| \leq \alpha\varepsilon_i\},
	\]
	where $\alpha>0, \varepsilon_i>0, f_i\in C(X,\mathbb{R}), \ i=1,\cdots,k$ and $\norm{f_i}\leq 1$ for each $i$, where $\norm{f_i}=\sup_{x\in X}|f_{i}(x)|$.
\end{definition}
The $f$-neighborhoods form a neighborhood base for the weak$^*$ topology on $M(X)$, which is the topology we use.

Following ideas in the Appendix part of \cite{ren2020topological}, we just need to show the proposition below.

\begin{proposition}\label{closedinvariant}
	Let $(X,G)$ be a dynamical system and $\{K_n\}$ a tempered F{\o}lner sequence and $\mu\in M(X,G)$. Suppose the system has the $g$-almost product property and $\mu$ verifies the conclusions of Proposition \ref{sepes1}. Let $0< h^{\prime} < h_{\mu}(X,G).$ Then there exists $\varepsilon^{\prime}>0$, such that for any neighborhood $C$ of $\mu$, there exists a $G$-invariant closed subset $Y\subset X$ satisfying the following properties.

	\begin{enumerate}
		\item There exists $n_{C}^{\prime} \in \mathbb{N},$ such that $\mathcal{E}_{K_n}(y) \in C$ for all $y\in Y$ and $n \geq n_{C}^{\prime}$;
		\item There exists $n_{C}^{\prime\prime} \in \mathbb{N},$ such that there exists a subset $\Gamma_n$ of $Y$ which is $(K_n,\varepsilon^{\prime})$-separated and $|\Gamma_n| \geq e^{|K_n|h^{\prime}}$ for all $n \geq n_{C}^{\prime\prime}.$
	\end{enumerate}
	In particular, $h_{top}(Y,G)\geq h^{\prime}$.
\end{proposition}

\begin{proof}
	Take $h^{\prime} < h^* < h_{\mu}(X,G)$. Given the neighborhood $C$ of $\mu$, take an $f$-neighborhood $F^{(1)}\subset C$ of $\mu$ with fixed $\{f_j,\varepsilon_j : j=1,2,\cdots,p\}.$ Denote $\varepsilon_{min} = \min\{\varepsilon_j \mid j=1,\cdots,p\}.$ Let $\delta^*, \varepsilon^*$ and $n^{*}_{F^{(1/5)}}$ correspond to $h^*$ in the conclusion of Proposition \ref{sepes1}. Set $n^* = n^{*}_{F^{(1/5)}}$. 
	
	Because $\{f_{j} \in C(X,\mathbb{R})\}$ are uniformly continuous on $X$, there exists $\triangle>0$ such that
	\[
	\triangle < \varepsilon^*/3 \text{ and } \rho(x,y) < \triangle \Longrightarrow \big|f_{j}(x) - f_{j}(y)\big| <\varepsilon_j/5
	\]
	for each $f_{j}$ associated with $F^{(1)}.$
	
Let $\{\gamma_n\}$ and $\{\tau_n\}$ be two sequence of decreasing positive numbers with $\gamma_n,\tau_n \rightarrow 0$ and $\tau_1 < \triangle.$  Let $g: (0,1) \rightarrow (0,1)$ be the mistake-density function as in the $g$-almost product property and $m:(0,1) \rightarrow F(G)\times (0,1)$ be the map in Definition \ref{almostspecification}. By Lemma \ref{tiling1} and Lemma \ref{newtiling}, there exist $N_1<n_{1,1,}<n_{1,t_1}<N_2<n_{2,1}<\cdots <n_{2,t_2}<\cdots$ and a congruent tilings $\{\mathcal{T}_k\}$ with shape sets $\{\mathcal{S}_k\}$ such that the following hold: Any $S\in \mathcal{S}_k$ can be $\gamma_k$-quasi tiled by $\{K_{n_1,1,},\cdots,K_{n_{k,t_k}}\}$ with tiling centers $\{C_{k,S,1},\cdots, C_{k,S,t_k}\}$. Denote the $\gamma_k$-quasi tile by
\begin{align}
\mathcal{T}_{k,S} = \{K_{{k,j}}c_{k,S,j}\mid c_{k,S,j}\in C_{k,S,j}, 1\leq j \leq t_k \}. \label{qusitile1}
\end{align}
Here we also assume $N_1 \geq n^{*}$. 

From Remark \ref{remark1}, we can modify $\mathcal{T}_{k,S}$ to get a new $\gamma_k$-quasi tile  
\begin{align}
\tilde{\mathcal{T}}_{k,S} = \{\tilde{K}_{c_{k,S,j}}c_{k,S,j}\mid K_{n_k,j}c_{k,S,j} \in \mathcal{T}_{k,S}\} \label{qusitile3}\end{align}
such that the following hold:
\begin{enumerate}
	\item Elements in $\tilde{\mathcal{T}}_{k,S}$ are pairwise disjoint and $\frac{|\tilde{K}_{c_{k,S,j}}|}{|K_{k,j}|}> 1- \gamma_k$;
	\item $\bigcup\tilde{\mathcal{T}}_{k,S} \subset S$ and $|\bigcup \tilde{\mathcal{T}}_{k,S}|>(1-\gamma_k)^2|S|$.
\end{enumerate}
Note that the sequence  obtained from  $\tilde{\mathcal{T}}_{k,S}$ is still a F{\o}lner sequence. By taking subsequence, we can assume for each $n\in\mathbb{N}$, any 
\[
A \in \bigcup_{k\geq n}\left\{ \tilde{K}_{c_{k,S,j}} \mid c_{k,S,j}\in C_{k,S,j}\ 1\leq j\leq t_k, S \in \mathcal{S}_k\right\}
\]

is $m(\tau_n)$-variant.

From the conclusion of Proposition \ref{sepes1}, we get the existence of  $(\delta^*,K_{c_{k,S,j}},\varepsilon^*)$-separated set $\Gamma_{c_{k,S,j}} \subset X_{K_{c_{k,S,j}},F^{(1/5)}}$ with $|\Gamma_{c_{k,S,j}}| \geq e^{|K_{c_{k,S,j}}|h^*}$ for $K_{c_{k,S,j}} \in \mathcal{T}_{k,S}$,  By Lemma \ref{eslemma1}, the set $\Gamma_{c_{k,S,j}}$ is also $(\delta^*/2,K_{c_{k,S,j}},\varepsilon^*)$-separated and contained in $X_{\tilde{K}_{c_{k,S,j}},F^{(2/5)}}.$

We mention that for $x_{c_k,S,j}\in \Gamma_{{c_{k,S,j}}} \subset X_{\tilde{K}_{c_{k,S,j}},F^{(2/5)}}$,
\begin{align}
|<f_j,\mathcal{E}_{\tilde{K}_{c_{k,S,j}}}(x_{c_{k,S,j}})> - <f_j,\mu>| < \frac{2}{5}\varepsilon_j .\label{appencal1}
\end{align}

For $Sd \in \mathcal{T}_k$, let
\begin{align*}
\Gamma(Sd) = \Gamma(S) & := \prod_{\tilde{K}_{c_{k,S,j}}c_{k,S,j}\in \tilde{\mathcal{T}}_{k,S}}\Gamma_{{c_{k,S,j}}}\\
& := \{\vec{x} = (x_{c_{k,S,j}}) \mid x_{c_{k,S,j}} \in \Gamma_{{c_{k,S,j}}}\}.
\end{align*}

Consider $Z_{F^{(1)},k}^{\#}$ defined by the requirement that $x\in Z_{F^{(1)},k}^{\#}$ if and only if for all $Sd\in \mathcal{T}_k$ there exists $\vec{x}\in \Gamma(Sd)$ such that 
\begin{align}
\rho_{\tilde{K}_{c_{k,S,j}}}(c_{k,S,j}dx,x_{c_{k,S,j}}) \leq \tau_k. \label{appenes1}
\end{align}
Since the system has the $g$-almost product property, the set $Z_{F^{1},k}^{\#}$ is not empty.

Let $\beta>0$ with $\beta < \frac{\varepsilon_{min}}{20}$ and $\beta< \frac{h^* - h^{\prime}}{3h^*}$. Choose $m_k$ large such that for $m\geq m_k$, $K_m$ is $(\cup \mathcal{S}_k,\frac{\beta}{|\cup\mathcal{S}_k|})$-invariant.

Let $m\geq m_k$. Define
\[
Y_{m,k} := \{x\in X\mid sx \in X_{K_m,F^{(4/5)}}, \forall s\in G\}.
\]
By the definition, $Y_{m,k}$ is a closed $G$-invariant subset. Next we will show that $Z_{F^{(1)},k}^{\#}\subset Y_{m,k}$. 

Take $s\in G$ and let $\Lambda_{K_{m}s} =\{T\in\mathcal{T}_k\mid T\subset K_{m}s\}$ and $\widetilde{K_{m}s}=\bigcup \Lambda_{K_{m}s}.$ 
Let 
\begin{align*}
I_n &=\{s\in F_n : \exists \ T\in \mathcal{T}_k \text{ such that } s\in T, T\cap (G\setminus F_n) \neq \emptyset \}\\
&\subset \bigcup\left\{Sc : \exists\  S\in \mathcal{S}_k, c\in G \text{ such that }Sc\cap F_n \neq \emptyset, Sc\cap (G\setminus F_n) \neq \emptyset\right\}\\
&\subset \bigcup \left\{(\cup \mathcal{S}_k)c : c\in \partial_{\cup \mathcal{S}_k} F_n\right\}.
\end{align*}
Hence $|I_n| \leq |\cup \mathcal{S}_k||\partial_{\cup \mathcal{S}_k} F_n| \leq \beta |F_n|$. Then we have 

\begin{align}
|\Lambda_{K_{m}s}|\geq (1-\beta)|F_n|.\label{appes1}
\end{align}

For $x\in Z_{F^{(1)},k}^{\#}$ and $\tilde{K}_{c_{k,S,j}}c_{k,S,j} \in \tilde{\mathcal{T}}_{k,S}$, by (\ref{appenes1}) and the choice of metric (\ref{metric}), we have
\begin{align}
D(\mathcal{E}_{\tilde{K}_{c_{k,S,j}}}(c_{k,S,j}dx), \mathcal{E}_{\tilde{K}_{c_{k,S,j}}}(x_{c_{k,S,j}})) \leq \tau_k. \label{appencal2}
\end{align}

For $Sd\in \Lambda_{K_{m}s}$, let $\tilde{S} = \bigcup\tilde{\mathcal{T}}_{k,S}.$  From (\ref{appencal1}) and (\ref{appencal2})
\begin{align}
\big|\frac{1}{|Sd|}\sum_{s\in Sd}f_{j}(sx) - \int f_j \dd\mu\big| & \leq \big|\frac{1}{|Sd|}\sum_{s\in Sd}f_{j}(sx) - \frac{1}{|Sd|}\sum_{s\in \tilde{S}d}f_{j}(sx)\big| \notag\\
&+ \big|\frac{1}{|Sd|}\sum_{s\in \tilde{S}d}f_{j}(sx) - \int f_j \dd\mu\big|\notag\\
&\leq \norm{f_j} \cdot 2\gamma_k + \frac{2}{5}\varepsilon_j + \norm{f_j}\tau_k \leq \frac{3}{5}\varepsilon_j.\label{appes2}
\end{align}
By (\ref{appes1}) and (\ref{appes2}),
\begin{align}
\big|\frac{1}{|K_{m}s|}\sum_{t\in K_{m}s}f_{j}(sx) - \int f_j \dd\mu\big| &\leq \big|\frac{1}{|K_{m}s|}\sum_{t\in K_{m}s}f_{j}(tx) - \frac{1}{|\tilde{K_{m}s}|}\sum_{t\in \widetilde{K_{m}s}}f_{j}(tx)\big| \notag\\
&+ \big|\frac{1}{|\widetilde{K_{m}s}|}\sum_{t\in \tilde{K_{m}s}}f_{j}(tx) - \int f_j \dd\mu\big|\notag \\
&\leq 2\beta \norm{f_j} + \frac{3}{5}\varepsilon_j \leq \frac{4}{5}\varepsilon_j. \label{appenre1}
\end{align}
From (\ref{appenre1}), we get $Z_{F^{(1)},k}^{\#}\subset Y_{m,k}$.

Define 
\[
Y := \cap_{m\geq m_k}Y_{m,k}.
\]
Then $Y$ is a non-empty closed $G$-invariant subset of $X.$ Set $n_{C}^{\prime}=m_k.$ For  $m\geq n_{C}^{\prime},$ we have $Y\subset Y_{m,k}$, which implies that for $y\in Y$, the measure $\mathcal{E}_{K_m}(y)\in F^{(4/5)} \subset C$. Then statement (1) is true.

Now we prove the statement(2) of this proposition. We set $n_{C}^{\prime\prime} = m_k$ and $\varepsilon^{\prime} = \frac{\varepsilon^*}{3}$.  Let $\Lambda_{K_{n}} =\{T\in\mathcal{T}_k\mid T\subset K_{n}\}$ and $\widetilde{K_{n}}=\cup \Lambda_{K_{n}}.$ From Lemma \ref{appes1},
\begin{equation}\label{apes7.1}
|\widetilde{K_{n}}|>(1-2\beta)|K_n|. 
\end{equation} 
For $S\in\mathcal{S}_k$, denote $\Gamma(S)=\prod_{i=1}^{t_k}\prod_{c_{k,S,i}\in C_{k,S,i}}\Gamma_{c_{k,S,i}}$. Denote $\Gamma(K_n) = \prod_{Sd\in \Lambda_{n}}\Gamma(S)$, where $\Lambda_n$ is as described in Lemma \ref{standardes1}. 

For each $n\geq n_{C}^{\prime\prime},$ we will consider a subset $Z_{n}^{\#} \subset Z_{F^{(1)},k}^{\#}$ with the following property: for each $\vec{x}=\{x_{Sd,c_{k,S,i}}\}\in \Gamma(K_n)$, there exists exact one point $x\in Z_{n}^{\#}$ such that
\begin{align}
\rho_{K_{c_{k,S,j}}}(c_{k,S,j}dx,x_{Sd,c_{k,S,j}}) \leq \tau_k. \label{map}
\end{align}
Define a map $\varPhi$ from $\Gamma(K_n)$ to $Z_{n}^{\#}$ such that $\varPhi(\vec{x})$ satisfies (\ref{map}). For $\vec{x} \neq \vec{y} \in \Gamma(K_n)$, we have 
\[
\rho_{K_n}(\varPhi(\vec{x}),\varPhi(\vec{y})) \geq \varepsilon^{*}-2\tau_k>\frac{\varepsilon^*}{3} = \varepsilon^{\prime}.
\]
Then $Z_{n}^{\#}$ is $(K_n,\varepsilon^{\prime})$-separated.
By the definition of $\Gamma(K_n)$ and $|\Gamma_{c_{k,S,i}}|\geq e^{h^{*}|K_{c_{k,S,i}}|}$
\begin{align}
|\Gamma(K_n)| &= \prod_{Sd\in \Lambda_{K_{n}}}\prod_{i=1}^{t_k}\prod_{c_{k,S,i}\in C_{k,S,i}}|\Gamma_{c_{k,S,i}}|\notag\\
&\geq e^{h^{*}\sum_{Sd\in\Lambda_{n}}\sum_{i=1}^{t_k}\sum_{c_{k,S,i}\in C_{k,S,i}}|K_{c_{k,S,i}}|} \notag\\
&\geq e^{h^{*}(1-2\beta)(1-4\gamma_k)|K_n|} \geq e^{h^{\prime}|K_n|}.
\end{align}
Thus statement (2) is true which implies $h_{top}(Y,G)\geq h^{\prime}$.

\end{proof}

{\bf Acknowledgements:} X. Ren is supported by National Natural Science Foundation of China(No. 11801261). W. Zhang is supported by National Science Foundation of Chongqing, China, with grant No. cstc2021jcyj-msxmX1045. Y. Zhang is supported by National Natural Science Foundation of China (No. 11871262), and Hubei Key Laboratory of Engineering Modeling and Scientific Computing in HUST. The authors would like to thank Prof. Xueting Tian and Prof. Yunhua Zhou for the helpful suggestions.
\bibliography{saturated_sets}
\end{document}